\documentclass[a4paper]{article}
\usepackage{latexsym,bm}
\usepackage{amssymb,amsmath,amsthm}
\usepackage[english,german]{babel}
\usepackage[colorlinks=true]{hyperref}
\usepackage{color}
\allowdisplaybreaks[4]

\newtheorem{theorem}{Theorem}[section]
\newtheorem{lemma}[theorem]{Lemma}
\newtheorem{remark}[theorem]{Remark}
\newtheorem{proposition}[theorem]{Propositon}
\newtheorem{definition}[theorem]{Definition}

\numberwithin{equation}{section}
\hypersetup{linkcolor=blue,urlcolor=red,citecolor=red}
\title{ Frequency-domain criterion on the stabilizability for infinite-dimensional linear control systems\thanks{This work was partially supported by the National Natural Science Foundation of
China under grant 11971022,12171359.}$\phantom{\;}^,$\thanks{Email: \;karl.kunisch@uni-graz.at (K. Kunisch),\; wanggs@yeah.net (G. Wang), \;huaiqiangyu@tju.edu.cn (H. Yu).}}
\author{\renewcommand{\thefootnote}{\arabic{footnote}}  Karl Kunisch\footnotemark[1]$\phantom{\;}^,$\footnotemark[2]\and \renewcommand{\thefootnote}{\arabic{footnote}}  Gengsheng Wang\footnotemark[3] \and \renewcommand{\thefootnote}{\arabic{footnote}} Huaiqiang Yu\footnotemark[1]$\phantom{\;}^,$\footnotemark[4]}

\date{}
\begin{document}
\selectlanguage{english}
\maketitle
\renewcommand{\thefootnote}{\arabic{footnote}} \footnotetext[1]{Johann Radon Institute for Computational and Applied Mathematics,  \"{O}AW, Altenbergerstrasse 69,
4040 Linz, Austria.}
    \footnotetext[2]{Institute of Mathematics and Scientific Computing,  University of Graz, Heinrichstrasse 36, 8010 Graz, Austria.}
\footnotetext[3]{Center for Applied Mathematics, Tianjin University, Tianjin, 300072, China.}
\footnotetext[4]{School of Mathematics, Tianjin University, Tianjin, 300354, China.}

\begin{abstract}
   A quantitative frequency-domain condition related to the
   exponential stabilizability for  infinite-dimensional linear control systems is presented.
   It is proven that this condition is  necessary and sufficient for  the stabilizability
   of special systems, while it is a necessary condition for the stabilizability in general.
   Applications are provided.
\end{abstract}


\section{Introduction}\label{yu-section-1}

Stabilization for linear control systems is one of the most important directions of control theory. How to determine whether a linear control system is  stabilizable is one of the largest concerns in this direction. Over the past half-century, researchers have obtained many useful results on this issue (see, for instance, \cite{Badra-Mitra-Raymond, Badra-Takahashi, Bardos-Lebeau-Rauch, Benaddi-Rao, Chen-Fulling-Sun, Fattorini, Huang-Wang, Liu-Wang-Xu-Yu, Liu, Liu-Liu-Rao, Ma-Wang-Yu, Raymond, Rebarber-Zwart-1998, Trelat-Wang-Xu, Weiss-Rebarber} and the references therein). These works have laid a solid foundation for the study of the stabilization of linear control systems. For finite-dimensional linear control systems, there is a well-known frequency-domain criterion
to determine the stabilizability, which is the Hautus test condition (see \cite{Hautus}). Unfortunately, this criterion may not be valid in infinite-dimensional settings. For infinite-dimensional linear control systems, researchers have been trying to obtain the corresponding frequency-domain criteria for  stabilization. This paper intends to provide a frequency-domain condition that is a necessary and sufficient condition for the stabilizability for special
linear control systems, while it is a necessary condition in general.
We begin by introducing the frequently used notation in this paper.

\subsection{Notation}
    Let $\mathbb{N}^+:=\mathbb{N}\setminus\{0\}$, $\mathbb{R}^+:=(0,+\infty)$ and $\mathbb{R}^-:=(-\infty,0)$.
    If  $\gamma\in\mathbb{R}$, we write
     $\mathbb{C}_{\gamma}^+:=\{z\in\mathbb{C}:\mbox{Re}\,z>\gamma\}$ and $\mathbb{C}^{-}_{\gamma}:=\{z\in\mathbb{C}:\mbox{Re}\,z<\gamma\}$. Let $\mathrm{i}$ be the unit imaginary number. If $S$ is a subset of $\mathbb{C}$, we denote its closure by $\overline{S}$. If $T>0$, we let $[T]:=\max\{n\in\mathbb{N}:n\leq T\}$. If $X$ is a Banach space, we denote its norm and dual space by $\|\cdot\|_X$ and $X^*$, respectively. If $X$ is a Hilbert space, we use $\langle\cdot,\cdot\rangle_X$ to denote its inner product. For Banach spaces $X_1$ and $X_2$, $\mathcal{L}(X_1;X_2)$
    denotes the space of all bounded linear operators from $X_1$ to $X_2$. We write $\mathcal{L}(X_1):=\mathcal{L}(X_1;X_2)$ if $X_1=X_2$.
    Given an unbounded (or bounded) linear operator $L$ from $X_1$ to $X_2$, its domain,
     kernel, adjoint operator, resolvent set, and spectrum are
           $D(L):=\{f\in X_1:Lf\in X_2\}$, $\mbox{Ker}(L):=\{f\in D(L):Lf=0\}$, $L^*$, $\rho(L)$ and $\sigma(L)$, respectively. Given two sets $\Lambda_1$ and $\Lambda_2$ in $X$, we let $\mbox{Span}\{\Lambda_1,\Lambda_2\}$ be the space spanned by the elements of $\Lambda_1$ and $\Lambda_2$. We use $C(\cdots)$ or $D(\cdots)$ to denote constants that depend on what is enclosed in the brackets.

\subsection{Control problem}
    Let $H$ and $U$ be two separable and complex Hilbert spaces. We consider the control system $[A,B]$, i.e.,
\begin{equation}\label{yu-5-4-1}
    y'(t)=Ay(t)+Bu(t),\;\;t\in\mathbb{R}^+,
\end{equation}
    (where $u\in L^2(\mathbb{R}^+;U)$) under the following assumptions:
\begin{enumerate}
     \item[$(A_1)$] Operator $A$ with its domain $D(A)$ generates  a $C_0$-semigroup $S(\cdot)$ on $H$.
\item[$(A_2)$] Operator $B$ belongs to $\mathcal{L}(U;H_{-1})$, where $H_{-1}$ is the completion of $H$ with respect to the norm $\|f\|_{-1}:=\|(\rho_0I-A)^{-1}f\|_{H}$ ($f\in H$), and
    $\rho_0\in \rho(A)\cap \mathbb{R}^+$ is arbitrarily fixed.

    \item[$(A_3)$] For each $T>0$, there is $C(T)>0$ such that
\begin{equation*}\label{yu-5-4-2-b}
    \int_0^T\|B^*S^*(t)\varphi\|_{U}^2dt\leq C(T)\|\varphi\|_{H}^2\;\;\mbox{for}\;\;\varphi\in D(A^*).
\end{equation*}
\end{enumerate}
\begin{remark}\label{yu-b-remark-5-11-1}
    There are several remarks on the above assumptions as follows:
\begin{enumerate}
  \item [(i)] In this paper, we write $H_1$ for the space $D(A^*)$ with the norm: $\|(\rho_0I-A^*)\varphi\|_H$, $\varphi\in D(A^*)$. (The space $H_{-1}$ is the dual space of $H_1$ with respect to the pivot space $H$, see \cite[ Section 2.9, Chapter 2]{Tucsnak-Weiss});
        $\widetilde{A}$ denotes the unique extension of  $A$ in
                 $\mathcal{L}(H;H_{-1})$, which is provided in the following manner
                 (see \cite[Proposition 2.10.3, Chapter 2]{Tucsnak-Weiss}):
\begin{equation}\label{yu-5-4-2}
    \langle \widetilde{A}\varphi,\psi\rangle_{H_{-1},H_1}=\langle \varphi, A^*\psi\rangle_{H}
    \;\;\mbox{for}\;\;\varphi\in H,\;\;\psi\in H_1;
\end{equation}
We let $\widetilde{S}(\cdot):=(\rho_0I-\widetilde{A})S(\cdot)(\rho_0I-\widetilde{A})^{-1}$,
which is the $C_0$-semigroup on $H_{-1}$ generated by $\widetilde{A}$ (with its domain $H$)
and an extension of $S(\cdot)$ (see \cite[Proposition 2.10.4, Chapter 2]{Tucsnak-Weiss}).

Moreover,  by assumption $(A_2)$, we have
 $B^*\in \mathcal{L}(H_1;U)$ and  $B^*(\rho_0I-A^*)^{-1}\in\mathcal{L}(H;U)$.

  \item [(ii)] By assumption $(A_3)$, we have that
  for any $u\in L^2(\mathbb{R}^+;U)$ and $y_0\in H$, system \eqref{yu-5-4-1}, which corresponds to this $u$ and with the initial condition $y(0)=y_0$, has a unique solution  in $C(\mathbb{R}^+;H)$,
  and this solution can be expressed by $y(t;y_0,u)=\widetilde{S}(t)y_0+\int_0^t\widetilde{S}(t-s)Bu(s)ds$,  $t\in\mathbb{R}^+$.
  (See \cite[Propositions 4.2.2, 4.2.5, Chapter 4]{Tucsnak-Weiss}.)
   \end{enumerate}
\end{remark}
   We are going to define the exponential/rapid stabilizability in $H$ for  system $[A,B]$,
   and `in $H$' will be omitted in what follows.

\begin{definition}\label{yu-def-5-5-2}
System $[A,B]$ is said to be
  exponentially stabilizable (or ``stabilizable" for short) if there is a constant $\alpha>0$,  a $C_0$-semigroup $\mathcal{S}_\dag(\cdot)$ on $H$ (with \color{red} its \color{black}
  generator $\mathcal{A}_\dag: D(\mathcal{A}_\dag)\subset H\to H$), and an operator
  $K_\dag\in \mathcal{L}(D(\mathcal{A}_\dag);U)$ such that
\begin{enumerate}
  \item [$(a)$] there is a constant $C_1>0$ such that $\|\mathcal{S}_\dag(t)\|_{\mathcal{L}(H)}
  \leq C_1e^{-\alpha t}$ for all  $t\in\mathbb{R}^+$;
\item[$(b)$] for each $x\in D(\mathcal{A}_\dag)$,  $\mathcal{A}_\dag x=\widetilde{A}x+B K_\dag x$, with $\widetilde{A}$ provided by \eqref{yu-5-4-2};
  \item [$(c)$] there is a constant  $C_2> 0$ such that $\|K_\dag\mathcal{S}_\dag(\cdot)x\|_{L^2(\mathbb{R}^+;U)}
  \leq C_2\|x\|_H$ for each $x\in D(\mathcal{A}_\dag)$.
\end{enumerate}
       $K_\dag$ and $\alpha$
       are  called the feedback law and a stabilizable decay rate (or ``a decay rate" for short), respectively.
  If $\alpha$,   $\mathcal{S}_\dag(\cdot)$, and
  $K_\dag$ exist,  system  $[A,B]$ is also said to be  stabilizable with decay rate $\alpha$.
\end{definition}

\begin{definition}\label{yu-def-5-5-3}
    System $[A,B]$ is
  rapidly (or completely) stabilizable if for each $\alpha>0$, $[A,B]$ is stabilizable with the  decay rate $\alpha$.
\end{definition}
\begin{remark}\label{yu-remark-7-8-1}
   Several notes on Definition \ref{yu-def-5-5-2} are provided as follows:
\begin{enumerate}
  \item [$(i)$] Definition \ref{yu-def-5-5-2} is originally from \cite{Flandoli-Lasiecka}, where the authors proved that
the solvability of the LQ problem: $V(y_0) = \inf_{u\in L^2(\mathbb{R}^+;U)}
\int_0^{+\infty}
 [\| y(t; u, y_0)\|^2_H+\| u(t)\|^2_U]dt$
(i.e., $V(y_0) < +\infty$  for all $y_0 \in  H$) implies  the  stabilizability of system
$[A,B]$ in the sense of Definition \ref{yu-def-5-5-2}. The reverse was  proven in  \cite[Proposition 3.9]{Liu-Wang-Xu-Yu}. Hence, the solvability of the above LQ problem is
equivalent to the  stabilizability of  system $[A,B]$ in the sense of
Definition \ref{yu-def-5-5-2}.

  \item [$(ii)$] If $B\in \mathcal{L}(U;H)$, the stabilizability of system $[A,B]$
   is defined as follows: there is a $K\in \mathcal{L}(H;U)$ such that $e^{(A+BK)t}$ is exponentially stable.
  Using the weak observability inequality in
  \cite{Liu-Wang-Xu-Yu, Trelat-Wang-Xu}, one can easily show that Definition \ref{yu-def-5-5-2}
  is an extension of the above definition if $B\in \mathcal{L}(U;H)$.

  \item [$(iii)$] If $[A,B]$ is stabilizable in the sense of Definition \ref{yu-def-5-5-2}, the feedback law
  can be constructed by the usual LQ theory.
\end{enumerate}

\end{remark}

\subsection{Motivation and novelty}
{\bf Motivation.} The stabilizability of system $[A,B]$ is equivalent to the
weak observability inequality. This equivalence can be summarized as the following lemma:

\begin{lemma}\label{yu-bb-remark-6-21-1}
\begin{enumerate}
  \item [$(i)$] The following statements are equivalent:
  \begin{enumerate}
    \item [$(a)$] System $[A,B]$ is  stabilizable.
    \item [$(b)$] There are constants  $\alpha>0$ and $C(\alpha)>0$ such that
\begin{equation}\label{yu-6-21-31}
    \|S^*(t)\varphi\|_H^2\leq C(\alpha)\left(\int_0^t\|B^*S^*(s)\varphi\|_U^2ds
    +e^{-\alpha t}\|\varphi\|_H^2\right)\;\;\mbox{for}\;\;t\in\mathbb{R}^+,\;\varphi\in H_1.
\end{equation}
    \item [$(c)$] There are constants $T>0$, $\delta\in (0,1)$ and $C>0$ such that
\begin{equation}\label{1.6-7-21w}
    \|S^*(T)\varphi\|_H^2\leq C\int_0^T\|B^*S^*(s)\varphi\|_U^2ds
    +\delta\|\varphi\|_H^2\;\;\mbox{for}\;\;\varphi\in H_1.
\end{equation}
  \end{enumerate}
  \item [$(ii)$] System $[A,B]$ is rapidly  stabilizable if and only if for each $\alpha>0$, there is $C(\alpha)>0$ that satisfies \eqref{yu-6-21-31}.
\end{enumerate}
\end{lemma}
\begin{remark}
 The equivalence of $(a)$ and $(c)$ above was originally obtained in \cite[Theorem 1]{Trelat-Wang-Xu} for the case that $B$ is bounded, and then extended to the case where $B$ is  admissible. The equivalence between $(b)$ and $(c)$ above was proved in
\cite[Proposition 3]{Ma-Wang-Yu}. The conclusion $(ii)$ of Lemma \ref{yu-bb-remark-6-21-1} was obtained
  in \cite[Theorem 1.1, Theorem 3.4]{Liu-Wang-Xu-Yu}.
\end{remark}
      Inequality \eqref{yu-6-21-31} (or \eqref{1.6-7-21w})
 can be considered a time-domain criterion for the stabilizability of system $[A,B]$.
Naturally, we would expect frequency-domain criteria  for system $[A,B]$.
For finite-dimensional
settings,
 the following criterion on the stabilizability
 is well known:
\begin{lemma}\label{yu-lemma-7-11-1} (\cite[Theorems 3, 4]{Hautus})
    Let  $H:=\mathbb{C}^{n}$, $U:=\mathbb{C}^m$, $A\in\mathbb{C}^{n\times n}$ and $B\in\mathbb{C}^{n\times m}$
    ($n,m\in\mathbb{N}^+$). Then, system $[A,B]$ is  stabilizable if and only if the pair $(A^*,B^*)$ satisfies  the condition
$$
    \mbox{(HSF)}:\;\;\mbox{Ker}(\lambda I-A^*,B^*)^{\top}=\{0\}\;\;\mbox{for} \;\; \lambda\in \sigma(A^*)\cap\overline{\mathbb{C}_0^+}.
$$
\end{lemma}
  The time-domain criterion and frequency-domain criterion have their own merits.
 Unfortunately, for general infinite-dimensional linear control systems,
 the equivalence in Lemma \ref{yu-lemma-7-11-1} may not be true. Counterexamples can be founded in \cite{Huang}
 (see also \cite[3.4 in Section 3, Chapter 4]{Engel-Nagel}). {\it Thus, it is natural to ask for frequency-domain criteria/conditions on the stabilization for system $[A,B]$ in our setting.}
\vskip 5pt
\noindent{\bf Novelty.} We provide a quantitative frequency-domain  condition
that can be considered as a generalization of condition (HSF) in Lemma \ref{yu-lemma-7-11-1}.
We prove that this condition  is a criterion on the stabilization for system $[A,B]$ under additional
conditions beyond $(A_1)$-$(A_3)$,
 while it  is
a necessary condition for the stabilization in general. Our method of proving these results uses the weak observability inequality in Lemma \ref{yu-bb-remark-6-21-1}. Combining time domain with
frequency domain concepts appears to be novel for
infinite dimensional stabilization problems.

\subsection{Main results}
   To state our main results, we need the following definition:
   \begin{definition}\label{yu-definition-6-9-1}
   $(i)$ The pair
     $(A^*,B^*)$ satisfies (HESI)$_\beta$ for $\beta>0$  if there exists $C(\beta)>0$ such that
\begin{equation}\label{yu-6-9-1}
    \|\varphi\|^2_H\leq \frac{C(\beta)}{(\mbox{Re}\,\lambda+\beta)^2}\left(\|(\lambda I-A^*)\varphi\|_H^2+\|B^*\varphi\|_U^2\right)
    \;\;\mbox{for}\;\;\lambda\in\mathbb{C}_{-\beta}^+,\;
    \varphi\in H_1,
\end{equation}
or equivalently, there exists $C(\beta)>0$ such that
\begin{equation}\label{yu-6-12-bb-0}
    \|\varphi\|_H^2\leq C(\beta)
    (\|(\lambda I-A^*)\varphi\|^2_H+\|B^*\varphi\|_U^2)\;\;\mbox{for}\;\;\lambda\in\mathbb{C}_{-\beta}^+,\;\varphi\in H_1.
\end{equation}

   $(ii)$  The pair $(A^*,B^*)$ satisfies (HESI) if it
         satisfies (HESI)$_\beta$ for each $\beta>0$.
\end{definition}
\begin{remark}
    Several  notes on Definition \ref{yu-definition-6-9-1} are as follows:
\begin{enumerate}
\item [$(i)$] `(HESI)' is the abbreviation of `\textbf{H}autus test condition for \textbf{E}xponential \textbf{S}tabilizability of \textbf{I}nfinite-dimensional systems'.
    The subscript $\beta$  is used to characterize the boundary of frequency-domain  that appears in (\ref{yu-6-9-1}) (or (\ref{yu-6-12-bb-0})).
\item [$(ii)$] The constants in (HESI)$_\beta$, \eqref{yu-6-9-1} and \eqref{yu-6-12-bb-0} are allowed to be different. The proof of the equivalence between \eqref{yu-6-9-1} and \eqref{yu-6-12-bb-0} will be given in Appendix of the paper (see Proposition \ref{yu-prop-6-3-1}).
    It is worth mentioning that \eqref{yu-6-9-1} is more sharp than \eqref{yu-6-12-bb-0} to describe the optimal decay rate of system $[A,B]$
(see Remark \ref{yu-remark-7-31-1}).

  \item [$(iii)$] (HESI)$_\beta$ (for some $\beta>0$)  is a type of quantitative frequency-domain condition.
  The connection to (HSF) (in Lemma \ref{yu-lemma-7-11-1}) is as follows: if $(A^*,B^*)$ satisfies (HESI)$_\beta$ for some $\beta>0$, then
    $\mbox{Ker}(\lambda I-A^*,B^*)^\top=\{0\}$ for each $\lambda\in\mathbb{C}_{-\beta}^+$.

  \item [$(iv)$]  (HESI)$_\beta$ (for some $\beta>0$) can be considered an
    extension of the classical frequency-domain criterion for the exponential stability of the $C_0$-semigroup
    $S(\cdot)$. Indeed, if $B^*=0$ and the pair $(A^*,0)$ satisfies (HESI)$_\beta$ for some $\beta>0$, then it follows from \eqref{yu-6-9-1} that
    $\mathbb{C}_{-\beta}^+\subset \rho(A^*)$ and
    $\sup_{\lambda\in\mathbb{C}_{-\alpha}^+}\|(\lambda I-A^*)^{-1}\|_{\mathcal{L}(H)}<+\infty$ for each $\alpha\in(0,\beta)$. Thus, by  \cite[Theorem 1.11, Chapter V]{Engel-Nagel}
    (see also \cite{Huang, Pruss}), we see that $S^*(\cdot)$ is exponentially stable, and so is $S(\cdot)$.
\end{enumerate}
\end{remark}

    The first two main theorems show that under additional conditions beyond $(A_1)$-$(A_3)$,
    (HESI)$_\beta$ for some $\beta>0$ is a necessary and sufficient condition on
    the stabilizability of system $[A,B]$.

\begin{theorem}\label{yu-theorem-6-30-1}
    Suppose that assumptions $(A_1)$-$(A_3)$ hold. Further assume that $A$ is a normal operator, and for each $\gamma>0$, $\sigma(A)\cap \mathbb{C}_{-\gamma}^+$ is bounded. Then, the following conclusions are true:
\begin{enumerate}
  \item [$(i)$] System $[A,B]$ is  stabilizable if and only if  the pair $(A^*,B^*)$ satisfies (HESI)$_\beta$ for some $\beta>0$.
  \item [$(ii)$] System $[A,B]$ is rapidly stabilizable if and only if the pair $(A^*,B^*)$ satisfies (HESI).
\end{enumerate}
\end{theorem}

\begin{remark}  \label{remark 1.10-7-24w}
Two notes on the assumptions in Theorem \ref{yu-theorem-6-30-1} are as follows:
\begin{enumerate}
  \item [$(i)$] The application of normal operators in partial differential equations covers a considerable area. It includes self-adjoint operators and linear partial differential operators on the entire space with (real/complex-valued) constant coefficients among others (see \cite[Theorem 13.24, Chapter 13]{Rudin}).

  \item [$(ii)$] If the semigroup $S(\cdot)$ is uniformly continuous/analytical/differentiable/compact, then $\sigma(A)\cap \mathbb{C}_{-\gamma}^+$ is bounded for each $\gamma>0$
      (see \cite[Sections 2.3-2.6, Chapter 2]{Pazy}).
\end{enumerate}
\end{remark}

\begin{theorem}\label{yu-theorem-5-20-1}
    Suppose that the assumptions $(A_1)$-$(A_3)$ hold. Assume that one of the following assumptions holds:
  \begin{enumerate}
    \item[$(a)$] The semigroup $S^*(\cdot)$ is uniformly bounded;
    \item[$(b)$] For each $\alpha>0$, there are two closed subspaces $Q_1:=Q_1(\alpha)$ and $Q_2:=Q_2(\alpha)$ of $H$ (depending on $\alpha$) such that $(b_1)$ $H=Q_1\oplus Q_2$;  $(b_2)$ $Q_1$ and $Q_2$ are  invariant subspaces of $S^*(\cdot)$;
    $(b_3)$ $A^*|_{Q_1}$, the restriction of $A^*$ on $Q_1$ is bounded and satisfies that  $\sigma(A^*|_{Q_1})\subset\mathbb{C}_{-\alpha}^+$; $(b_4)$  $S^*(\cdot)|_{Q_2}$,   the restriction of $S^*(\cdot)$ on $Q_2$ is  exponentially stable.
  \end{enumerate}
   Then, system $[A,B]$ is stabilizable if and only if the pair $(A^*,B^*)$ satisfies (HESI)$_\beta$ for some $\beta>0$.
\end{theorem}

\begin{remark}\label{yu-remark-6-9-1}
    Several notes on Theorem \ref{yu-theorem-5-20-1} are  as follows:
\begin{enumerate}
\item[$(i)$] If $A$ is skew-adjoint (i.e., $A=-A^*$), then the semigroup $S^*(\cdot)$ is uniformly bounded.
 \item[$(ii)$] If the semigroup $S(\cdot)$ is compact, then for each $\alpha>0$, we can find two closed subspaces $Q_1$ and $Q_2$ of $H$ such that $Q_1$ is finite-dimensional;  conditions $(b_1)$-$(b_4)$ in  $(b)$ of Theorem \ref{yu-theorem-5-20-1} hold. (See
    \cite[Section 3, Chapter IV, Section 3, ChapterV]{Engel-Nagel}.)
   \item[$(iii)$] In assumption $(b)$ of Theorem \ref{yu-theorem-5-20-1}, it is not required that
   $Q_1$ has a finite dimension.
 \item[$(iv)$] By the  spectral theorem (e.g., \cite[Theorem 6.47, Chapter 6]{Kubrusly}), one can check that if $A$ satisfies the assumptions of Theorem \ref{yu-theorem-6-30-1}, then it also satisfies assumption $(b)$ in
      Theorem \ref{yu-theorem-5-20-1}. Thus, $(i)$ in Theorem \ref{yu-theorem-6-30-1} can be regraded as an application of Theorem \ref{yu-theorem-5-20-1}.  But, there are some difficulties in applying Theorem \ref{yu-theorem-5-20-1} to characterize the rapid stabilizability. Indeed, when $A$ satisfies the assumptions of Theorem \ref{yu-theorem-6-30-1}, $A+\beta I$ also satisfies these assumptions for each $\beta>0$. This property plays a crucial role in the proof of $(ii)$ of Theorem \ref{yu-theorem-6-30-1}. It is not known  to us
      if assumption $(b)$ in Theorem \ref{yu-theorem-5-20-1} can guarantee this property for the general cases.
\end{enumerate}
\end{remark}

The last main result of this paper shows that  (HESI)$_\beta$ (for some $\beta>0$)/(HESI) is a necessary condition on the exponentially/rapidly stabilizability for system $[A,B]$, with $(A_1)$-$(A_3)$ holding.
\begin{theorem}\label{yu-theorem-7-3-1}
    Under assumptions $(A_1)$-$(A_3)$, the following conclusions are true:
\begin{enumerate}
  \item [$(i)$] If system $[A,B]$ is  stabilizable, then
   the pair $(A^*,B^*)$ satisfies
    (HESI)$_\beta$ for some $\beta>0$.
  \item [$(ii)$] If system $[A,B]$ is  rapidly stabilizable, then
   the pair $(A^*,B^*)$ satisfies
    (HESI).
\end{enumerate}

\end{theorem}

\subsection{Related works}\label{yu-sec-1.3}

\begin{itemize}
    \item H. O. Fattorini in \cite{Fattorini} established a frequency-domain condition similar to Lemma \ref{yu-lemma-7-11-1}
    for the special infinite-dimensional setting, where $B\in\mathcal{L}(U;H)$; system $[A,B]$
     can be decomposed  into two decoupled subsystems: one is in a finite-dimensional subspace and controllable; the other is in an infinite-dimensional subspace and exponentially stable.
     For more studies in this direction, we refer the readers to
     \cite{Badra-Takahashi, Curtain-Zwart, Raymond} and the references therein.
     The setting in \cite{Fattorini} is covered by our setting, where $(A_1)$-$(A_3)$ and
       $(b)$ in Theorem \ref{yu-theorem-5-20-1} hold (see $(ii)$ in Remark \ref{yu-remark-6-9-1}).
                 Thus,  Theorem \ref{yu-theorem-5-20-1}
      can be considered an extension of the related result in \cite{Fattorini}.

\item K. Liu proved in \cite{Liu} that if $A$ is skew-adjoint
 and
$B\in\mathcal{L}(U;H)$,  then  system $[A,B]$ is stabilizable if and only if
 $\mathrm{i}\mathbb{R}\in \rho(A-BB^*)$ and
$\sup_{\omega\in\mathbb{R}}\|(\mathrm{i}\omega I-A+BB^*)^{-1}\|_{\mathcal{L}(H)}<+\infty$. Later,
 Q. Zhou and M. Yamamoto in \cite{Zhou-Yamamoto} obtained that  if $A$ is skew-adjoint
 and
$B\in\mathcal{L}(U;H)$, then system $[A,B]$ is stabilizable if and only if there was $C>0$ such that
$$
    \|\varphi\|^2_H\leq C(\|(\mathrm{i}\omega I-A^*)\varphi\|^2_H+\|B^*\varphi\|^2_U)
    \;\;\mbox{for}\;\;\omega\in\mathbb{R},\;\;\varphi\in H_1.
$$
These are  two  frequency-domain criteria on the stabilizability for the case $A=-A^*$
 and
$B\in\mathcal{L}(U;H)$.
However, the condition that $A=-A^*$ implies condition $(a)$ in Theorem \ref{yu-theorem-5-20-1}.
Thus, Theorem \ref{yu-theorem-5-20-1} can also be considered an extension of the above frequency-domain criteria.

\item R. Rebarber and H. Zwart in \cite{Rebarber-Zwart-1998} introduced the  concept of the open-loop stabilizability for system $[A,B]$ and provided necessary conditions in the frequency-domain
    for the  open-loop stabilizability in infinite dimensional settings. Such
    stabilizability is defined as follows:  \emph{If there is $\sigma>0$ such that for each $y_0\in H$, there is a control $u\in \mathcal{D}'(\mathbb{R}^+;U)(=(C_0^\infty(\mathbb{R}^+;U))')$ such that the
    solution (in the sense of distribution)
    $y(\cdot;u,y_0)$ to system \eqref{yu-5-4-1} (with the initial condition $y(0)=y_0$ and the control $u$)  satisfies
    $e^{\sigma\cdot}y(\cdot;u,y_0)\in L^2(\mathbb{R}^+;H)$, then system $[A,B]$ is called open-loop stabilizable.} Clearly,  the open-loop stabilizability is weaker than the closed-loop stabilizability.
    (Several examples that are open-loop stabilizable but not closed-loop stabilizable were provided in
    \cite{Rebarber-Zwart-1998}.)
    For more studies on the open-loop stabilizability in infinite dimensional settings, we refer the readers to \cite{Rebarber-Zwart-1993, Weiss-Rebarber, Zwart-1989, Zwart-1989-1, Zwart-1990-1}.
    Our condition (HESI)$_\beta$ is partially inspired by \cite{Rebarber-Zwart-1998} and related  works.

\end{itemize}

\subsection{Plan of this paper}
    This paper is organized as follows: In Section 2, we present
    some criteria on the stabilizability for system
     $[A,B]$. In Section 3, we prove our main theorems. In Section  4, we provide selected applications.
Section 5 is Appendix.

\section{Other criteria on the stabilization}\label{yu-sec-2}
    This  section provides  criteria on   stabilizability for system $[A,B]$ from the perspective of  integral transformation.
    Although they are not easily verifiable, they play  important roles in the
     proofs of our main theorems.

     To present them, we must introduce the following function spaces:
       Let $X$ be a sparable, complex Hilbert space and $\alpha\geq 0$. For each
       open and connected subset  $S\subset \mathbb{C}$, we let $\mathcal{H}(S;X)$ be the set of all $X$-valued holomorphic functions on $S$.
    We define  the following Hardy space:
\begin{equation*}
    \mathcal{H}^2(\mathbb{C}^+_{-\alpha};X)
    :=\left\{f\in \mathcal{H}(\mathbb{C}^+_{-\alpha};X):\mathcal{M}_\alpha(f)<+\infty\right\}
    \;\;\mbox{with}\; \mathcal{M}_\alpha(f):=
    \sup_{\omega_1>-\alpha}\int_{\mathbb{R}}
    \|f(\omega_1+\textrm{i}\omega_2)\|^2_Xd\omega_2,
\end{equation*}
    and  the following weighted $L^2$-space:
\begin{equation*}\label{yu-3-17-1-b}
    L_\alpha^2(\mathbb{R}^+;X):=\{h:\mathbb{R}\to X:e^{\alpha\cdot}h(\cdot)\in L^2(\mathbb{R};X);
    \;\;h(\cdot)=0\;\;\mbox{on}\;\;\mathbb{R}^-\},
\end{equation*}
    with the inner product:
\begin{equation*}\label{yu-3-17-2-b}
    \langle g,h\rangle_{L^2_\alpha(\mathbb{R}^+;X)}:=\int_{\mathbb{R}^+} e^{2\alpha t}g(t)h(t)dt.
   \end{equation*}
    One can easily check that $L_\alpha^2(\mathbb{R}^+;X)$ is a Hilbert space and continuously embedded into
    $L^2(\mathbb{R}^+;X)$.
    Throughout this paper, we extend each $f\in L^2(\mathbb{R}^+;X)$ over $\mathbb{R}$ by setting it to be zero
    over $\mathbb{R}^-$ (we denote this extension in same way). Thus, we have  $L^2_0(\mathbb{R}^+;X)=L^2(\mathbb{R}^+;X)$.

    \begin{theorem}\label{yu-main-theorem-1}
    Suppose that $(A_1)$-$(A_3)$ hold. Then, the following statements are equivalent:
\begin{enumerate}
  \item [$(i)$] System $[A,B]$ is  stabilizable.
  \item [$(ii)$] There is $\alpha>0$ such that for each $\beta\in[0,\alpha)$ and  $y_0\in H$, there  is $(\xi(\cdot;y_0),\eta(\cdot;y_0))\in \mathcal{H}^2(\mathbb{C}^+_{-\beta};H)\times \mathcal{H}^2(\mathbb{C}^+_{-\beta};U)$ such that
\begin{equation}\label{yu-3-24-18}
    \langle \xi(\lambda;y_0),(\bar{\lambda}I-A^*)\varphi\rangle_H+\langle \eta(\lambda;y_0),B^*\varphi\rangle_U=\langle y_0,\varphi\rangle_H
    \;\;\mbox{for}\;\;\lambda\in\mathbb{C}_{-\beta}^+,\;\varphi\in H_1;
\end{equation}
   \begin{equation}\label{yu-3-24-18-b}
    \|\xi(\lambda;y_0)\|_H\leq \frac{C(\beta)}{\mbox{Re}\,\lambda+\beta}\|y_0\|_H,\;\;
    \|\eta(\lambda;y_0)\|_U\leq \frac{D(\beta)}{\sqrt{\mbox{Re}\,\lambda+\beta}}\|y_0\|_H \;\;\mbox{for}\;\;\lambda\in\mathbb{C}^+_{-\beta},
\end{equation}
     where $C(\beta)> 0$ and $D(\beta)> 0$ are two constants independent of $y_0$.
\item[$(iii)$] There is $\beta\geq 0$ such that for each $y_0\in H$, there  is $(\xi(\cdot;y_0),\eta(\cdot;y_0))\in \mathcal{H}^2(\mathbb{C}^+_{-\beta};H)\times \mathcal{H}^2(\mathbb{C}^+_{-\beta};U)$ satisfying \eqref{yu-3-24-18}.
\end{enumerate}
\end{theorem}

To prove Theorem \ref{yu-main-theorem-1}, we need the following lemmas:
 The first lemma contains several  results quoted from  \cite[Section A.6.3, Chapter A]{Curtain-Zwart}:
\begin{lemma}\label{yu-lemma-6-9-1}
    Let $\alpha\geq 0$. Then, the following statements are true:
\begin{enumerate}
  \item [$(i)$] $f\in \mathcal{H}^2(\mathbb{C}_{-\alpha}^+;X)$ if and only if there is a unique $h\in L^2_\alpha(\mathbb{R}^+;X)$ such that
      $f(\lambda)=\int_{\mathbb{R}^+}e^{- t\lambda}h(t)dt$ ($\lambda\in \mathbb{C}_{-\alpha}^+$), i.e., $\mathcal{H}^2(\mathbb{C}^+_{-\alpha};X)$ and $L_\alpha^2(\mathbb{R}^+;X)$ are linear isomorphic. Moreover, $\frac{1}{2\pi}\mathcal{M}_\alpha(f)=\int_{\mathbb{R}^+}e^{2\alpha t}\|h(t)\|_X^2dt$.
  \item [$(ii)$] For each $f\in \mathcal{H}^2(\mathbb{C}^+_{-\alpha};X)$, there is a unique $f_\alpha^*\in L^2(\mathbb{R};X)$ such that
    $\lim_{\omega_1\to -\alpha}\|f(\omega_1+\mathrm{i}\cdot)-f_\alpha^*\|_{L^2(\mathbb{R};X)}=0$.
      Moreover,
    $\|f_\alpha^*\|_{L^2(\mathbb{R};X)}^2=
    \mathcal{M}_\alpha(f)$.
  \item [$(iii)$] $\mathcal{H}^2(\mathbb{C}^+_{-\alpha};X)$, with inner product
    $\langle f,g\rangle_{\mathcal{H}^2(\mathbb{C}^+_{-\alpha};X)}:=\langle f_\alpha^*,g^*_\alpha\rangle_{L^2(\mathbb{R};X)}$ ($f,g\in \mathcal{H}^2(\mathbb{C}^+_{-\alpha};X)$),
    is a Hilbert space.
\end{enumerate}
\end{lemma}
 The second lemma is as follows, which is clear if $B\in\mathcal{L}(U;H)$. However,  for our framework, we do not find accurate literature that provides its proof. Thus, we provide it  for
   the completeness of the paper.
 \begin{lemma}\label{yu-lemma-5-11-1}
    Suppose that $(A_1)$-$(A_3)$ hold. If system $[A,B]$ is  stabilizable with decay rate $\alpha>0$, then
    for each $\beta\in(0,\alpha)$, system $[A+\beta I,B]$ is  stabilizable.
\end{lemma}
\begin{proof}
We arbitrarily fix $\beta\in (0,\alpha)$. Since  system $[A,B]$ is  stabilizable with decay rate $\alpha>0$, we can use  the same method as that used in \emph{Step 1} of the proof of \cite[Theorem 3.4]{Liu-Wang-Xu-Yu} to  find a positive constant $C(\alpha)$ such that
\begin{equation*}\label{yu-5-11-1}
    \|S^*(t)\varphi\|_H^2\leq C(\alpha)\Big(\int_0^t\|B^*S^*(s)\varphi\|^2_Uds
    +e^{-2\alpha t}\|\varphi\|^2_H\Big)\;\;\mbox{for}\;\;t\in\mathbb{R}^+,\;\varphi\in H_1.
\end{equation*}
   This  yields  that for each  $t\in\mathbb{R}^+$,
\begin{equation}\label{yu-5-11-2}
    \|S^*_\beta(t)\varphi\|_H^2\leq C(\alpha)e^{2\beta t}\int_0^t\|B^*S_{\beta}^*(s)\varphi\|_H^2ds
    +C(\alpha)e^{-2(\alpha-\beta)t}\|\varphi\|_H^2\;\;\mbox{for}\;\;\varphi\in H_1,
\end{equation}
where $S_\beta(\cdot)$ is the $C_0$-semigroup on $H$ generated by $A_\beta:=A+\beta I$ with its domain
    $D(A_\beta)=D(A)$. Let $T>0$ satisfy $\delta_T:=C(\alpha)e^{-2(\alpha-\beta)T}<1$. Then, it follows from \eqref{yu-5-11-2} that
\begin{equation*}\label{yu-5-11-3}
    \|S^*_\beta(T)\varphi\|_H^2\leq C(\alpha)e^{2\beta T}\int_0^T\|B^*S_{\beta}^*(s)\varphi\|_H^2ds
    +\delta_T\|\varphi\|_H^2\;\;\mbox{for}\;\;\varphi\in H_1.
\end{equation*}
  With  Lemma \ref{yu-bb-remark-6-21-1}, the above shows
  that system $[A_\beta,B]$ is  stabilizable. This completes the proof.
\end{proof}

\begin{proof}[Proof of Theorem \ref{yu-main-theorem-1}]
  We organize the proof into the following steps.
\vskip 5pt
  \noindent  \emph{Step 1. We prove $(i)\Rightarrow (ii)$.}

  Suppose that system $[A,B]$ is stabilizable with decay rate $\alpha>0$.
  We arbitrarily fix  $\beta\in [0,\alpha)$. According to  Lemma \ref{yu-lemma-5-11-1},  system $[A_\beta,B]$ (where $A_\beta:=A+\beta I$) is  stabilizable, i.e.,
      there is $\gamma>0$, a $C_0$-semigroup $\mathcal{S}_{\beta,\gamma}(\cdot)$ on $H$ with the generator $\mathcal{A}_{\beta,\gamma}: D(\mathcal{A}_{\beta,\gamma})\subset H\to H$, and a $K_{\beta,\gamma}\in \mathcal{L}(D(\mathcal{A}_{\beta,\gamma});U)$ such that
     \begin{enumerate}
       \item [(a1)]  there is  $C_{\beta,\gamma,1}> 0$ such that $\|\mathcal{S}_{\beta,\gamma}(t)\|_{\mathcal{L}(H)}\leq C_{\beta,\gamma,1}e^{-\gamma t}$ for all $t\in\mathbb{R}^+$;
       \item [(b1)]  for each $x\in D(\mathcal{A}_{\beta,\gamma})$, $\mathcal{A}_{\beta,\gamma}x=\widetilde{A}_\beta x+BK_{\beta,\gamma}x$, where
           $\widetilde{A}_\beta:=\widetilde{A}+\beta I$, with $\widetilde{A}$  provided by (\ref{yu-5-4-2});
       \item [(c1)]  there is  $C_{\beta,\gamma,2}> 0$ such that $\|K_{\beta,\gamma}\mathcal{S}_{\beta,\gamma}(\cdot)x\|_{L^2(\mathbb{R}^+;U)}\leq C_{\beta,\gamma,2}\|x\|_H$ for each $x\in D(\mathcal{A}_{\beta,\gamma})$.
     \end{enumerate}
    From these properties we deduce the following facts:
\begin{enumerate}
  \item [(O1)] Based on (b1) and (c1), for each $y_0\in D(\mathcal{A}_{\beta,\gamma})$,
\begin{equation}\label{w2.4-7-20}
  u_{K_{\beta,\gamma}}(\cdot;y_0):=K_{\beta,\gamma}\mathcal{S}_{\beta,\gamma}(\cdot)y_0
  =K_{\beta,\gamma}y_{K_{\beta,\gamma}}(\cdot;y_0)
  \in L^2_0(\mathbb{R}^+;U).
\end{equation}
    and
\begin{equation*}\label{yu-5-8-10}
    y_{K_{\beta,\gamma}}(t;y_0):=\mathcal{S}_{\beta,\gamma}(t)y_0
    =\widetilde{S}_\beta(t)y_0+\int_0^t\widetilde{S}_\beta(t-s)
    Bu_{K_{\beta,\gamma}}(s;y_0)ds
    \;\;\mbox{for all}\;\;t\in\mathbb{R}^+,
\end{equation*}
    where $\widetilde{S}_{\beta}(\cdot):=e^{\beta\cdot}\widetilde{S}(\cdot)$  is the $C_0$-semigroup generated by $\widetilde{A}_\beta$ on $H_{-1}$ (see $(i)$ in Remark \ref{yu-b-remark-5-11-1}).
  \item [(O2)] We arbitrarily fix $y_0\in D(\mathcal{A}_{\beta,\gamma})$.
  By (a1) and (\ref{w2.4-7-20}), we obtain that if  $t\in\mathbb{R}^+$,
\begin{eqnarray}\label{yu-5-9-2}
    \|u_{K_{\beta,\gamma}}(t;y_0)\|_U&=&\|K_{\beta,\gamma}y_{K_{\beta,\gamma}}(t;y_0)\|_U
    =\|K_{\beta,\gamma}\mathcal{S}_{\beta,\gamma}(t)y_0\|_U\nonumber\\
    &\leq&\|K_{\beta,\gamma}(\rho_1 I-\mathcal{A}_{\beta,\gamma})^{-1}\|_{\mathcal{L}(H;U)}
    \|\mathcal{S}_{\beta,\gamma}(t)(\rho_1 I-\mathcal{A}_{\beta,\gamma})y_0\|_H\nonumber\\
    &\leq& C_{\beta,\gamma,1}e^{-\gamma t}\|K_{\beta,\gamma}(\rho_1 I-\mathcal{A}_{\beta,\gamma})^{-1}\|_{\mathcal{L}(H;U)}\|(\rho_1 I-\mathcal{A}_{\beta,\gamma})y_0\|_H,
\end{eqnarray}
 where
     $\rho_1\in \rho(\mathcal{A}_{\beta,\gamma})\cap \mathbb{R}^+$, and we use that
     $K_{\beta,\gamma}\in\mathcal{L}(D(\mathcal{A}_{\beta,\gamma});U)$.
  \item [(O3)] We arbitrarily fix  $y_0\in D(\mathcal{A}_{\beta,\gamma})$ and define
 \begin{equation}\label{yu-5-10-3}
    \xi(\lambda;y_0):=\int_{\mathbb{R}^+}e^{-\lambda t}y_{K_{\beta,\gamma}}(t;y_0)dt;\;\;\eta(\lambda;y_0):=-\int_{\mathbb{R}^+}e^{-\lambda t}u_{K_{\beta,\gamma}}(t;y_0)dt,\;\;\lambda\in \mathbb{C}^+_{0}.
\end{equation}
 Considering the above two functions, we have the following conclusions: First,
 by \eqref{yu-5-10-3} and $(i)$ of Lemma \ref{yu-lemma-6-9-1}, we see that
  $ \xi(\cdot;y_0)\in\mathcal{H}^2(\mathbb{C}^+_{0};H)$ and $\eta(\cdot;y_0)\in\mathcal{H}^2(\mathbb{C}^+_{0};U)$.
  Second, based on \eqref{w2.4-7-20}-\eqref{yu-5-10-3}, (a1) and (c1),
  we find
\begin{equation}\label{yu-5-10-10}
    \|\xi(\lambda;y_0)\|_H\leq \frac{C_{\beta,\gamma,1}}{\mbox{Re}\,\lambda}\|y_0\|_H
    \;\;\mbox{and}\;\;
    \|\eta(\lambda;y_0)\|_U\leq \frac{C_{\beta,\gamma,2}}{\sqrt{\mbox{Re}\,\lambda}}\|y_0\|_H\;\;
    \mbox{for}\;\;\lambda\in\mathbb{C}^+_0.
\end{equation}
  Third, based on (a1), (c1), and  Lemma \ref{yu-lemma-6-9-1}, we obtain
\begin{equation}\label{yu-5-11-10}
    \|\xi(\cdot;y_0)\|^2_{\mathcal{H}^2(\mathbb{C}^+_0;H)}=2\pi \int_{\mathbb{R}^+}\|y_{K_{\beta,\gamma}}(t;y_0)\|_H^2dt\leq \gamma^{-1}\pi C^2_{\beta,\gamma,1}\|y_0\|_H^2;
\end{equation}
    \begin{equation*}\label{yu-5-11-11}
    \|\eta(\cdot;y_0)\|^2_{\mathcal{H}^2(\mathbb{C}^+_0;U)}=2\pi \int_{\mathbb{R}^+}\|u_{K_{\beta,\gamma}}(t;y_0)\|_U^2dt\leq 2\pi C^2_{\beta,\gamma,2}\|y_0\|_H^2.
\end{equation*}
    \end{enumerate}

   Next, we  show Claim One:  For each $y_0\in H$,  there is $(\xi(\cdot;y_0),\eta(\cdot;y_0))
    \in\mathcal{H}^2(\mathbb{C}_0^+;H)\times \mathcal{H}^2(\mathbb{C}_0^+;U)$ that satisfies (\ref{yu-5-10-10}) and
\begin{equation}\label{yu-5-12-1}
    \langle \xi(\lambda;y_0),(\bar{\lambda}I-A^*_\beta)\varphi\rangle_H+\langle \eta(\lambda;y_0),B^*\varphi\rangle_U=\langle y_0,\varphi\rangle_H
    \;\;\mbox{for}\;\;\lambda\in\mathbb{C}_{0}^+,\;\varphi\in H_1.
\end{equation}
    The proof of Claim One will be organized using two cases.
\vskip 5pt
  \noindent  \emph{Case 1. We consider that $y_0\in D(\mathcal{A}_{\beta,\gamma})$.}

   First, based on (O3), we have (\ref{yu-5-10-10}) for this case.
   We now show \eqref{yu-5-12-1} for this case. For this purpose, we arbitrarily fix
   $y_0\in D(\mathcal{A}_{\beta,\gamma})$ and $\varphi\in H_1=D(A^*_\beta)$.
      From (O1) and the main theorem in \cite{Ball} (see also \cite[Theorem on Page 259]{Pazy}), we obtain
\begin{equation}\label{yu-5-8-11}
\begin{cases}
    \frac{d}{dt}\langle y_{K_{\beta,\gamma}}(t;y_0),\varphi\rangle_H=\langle y_{K_{\beta,\gamma}}(t;y_0),A^*_\beta\varphi\rangle_H
    +\langle u_{K_{\beta,\gamma}}(t;y_0),B^*\varphi\rangle_{U},&t\in\mathbb{R}^+,\\
    y_{K_{\beta,\gamma}}(0)=y_0.
\end{cases}
\end{equation}
    The combination of this result, (a1), and (\ref{yu-5-9-2}) yield
\begin{equation}\label{yu-5-9-1}
    e^{-\lambda \cdot}\frac{d}{dt}\langle y_{K_{\beta,\gamma}}(\cdot;y_0),\varphi\rangle_H\in L^1(\mathbb{R}^+;\mathbb{C})\;\;\mbox{for each}
    \;\;\lambda\in\mathbb{C}^+_{0}.
\end{equation}
   Now, (a1), (\ref{yu-5-9-1}), and (\ref{yu-5-10-3}) lead to
\begin{equation*}\label{yu-5-10-4}
    \int_{\mathbb{R}^+}e^{-\lambda t}\frac{d}{dt}\langle y_{K_{\beta,\gamma}}(t;y_0),\varphi\rangle_Hdt
    =\langle\lambda\xi(\lambda;y_0),\varphi\rangle_H -\langle y_0,\varphi\rangle_H
    =\langle\xi(\lambda;y_0),\bar{\lambda}\varphi\rangle_H -\langle y_0,\varphi\rangle_H,\; \lambda\in \mathbb{C}^+_{0}.
\end{equation*}
   From the above and (\ref{yu-5-8-11}), one can directly obtain that
    $(\xi(\cdot;y_0),\eta(\cdot;y_0))$ (which is defined by (\ref{yu-5-10-3})) satisfies
    (\ref{yu-5-12-1}).
\vskip 5pt
 \noindent   \emph{ Case 2. We consider that $y_0\in H$.}

     According to the density of $D(\mathcal{A}_{\beta,\gamma})$ in $H$, there is a sequence $\{y_0^n\}_{n\in\mathbb{N}^+}\subset D(\mathcal{A}_{\beta,\gamma})$ such that
    $y_0^n\to y_0$ in $H$ as $n\to+\infty$. Thus, $\{y_0^n\}_{n\in\mathbb{N}^+}$ is a Cauchy sequence in $H$.
    Moreover, by (\ref{yu-5-10-3}) and (\ref{yu-5-11-10}), we have
\begin{equation*}\label{yu-5-12-2}
    \|\xi(\cdot;y_0^n)-\xi(\cdot;y_0^m)\|_{\mathcal{H}^2(\mathbb{C}_0^+;H)}
    =\|\xi(\cdot;y_0^n-y_0^m)\|_{\mathcal{H}^2(\mathbb{C}_0^+;H)}\leq \sqrt{\gamma^{-1}\pi}C_{\beta,\gamma,1}\|y_0^n-y_0^m\|_H\;\;\forall\; n,m\in\mathbb{N}^+.
\end{equation*}
    Hence, $\{\xi(\cdot;y_0^n)\}_{n\in\mathbb{N}^+}$ is
    a Cauchy sequence in $\mathcal{H}^2(\mathbb{C}_0^+;H)$.
    Then, according to $(iii)$ in Lemma \ref{yu-lemma-6-9-1},
    there is  $\xi^*(\cdot)\in \mathcal{H}^2(\mathbb{C}_0^+;H)$ such that
\begin{equation}\label{yu-5-12-3}
    \xi(\cdot;y_0^n)\to \xi^*(\cdot)\;\;\mbox{in}
    \;\;\mathcal{H}^2(\mathbb{C}_0^+;H)\;\;\mbox{as}\;\;n\to+\infty.
\end{equation}
   By $(i)$ in Lemma \ref{yu-lemma-6-9-1}, we can find $h^*(\cdot)$ in $L^2(\mathbb{R}^+;H)$ such that
    $\xi^*(\cdot)=\int_{\mathbb{R}^+}e^{-\cdot t}h^*(t)ds$.
     The combination of this result, \eqref{yu-5-10-3}, (\ref{yu-5-12-3}) and $(i)$ in Lemma \ref{yu-lemma-6-9-1} yield that
     if we write $h^n(t):=y_{K_{\beta,\gamma}}(t;y^n_0)$, $t\in \mathbb{R}^+$ ($n\in \mathbb{N}^+$), then we have
\begin{equation*}\label{yu-5-12-5}
    \int_{\mathbb{R}^+}\|h^n(t)-h^*(t)\|_H^2dt=
    \frac{1}{2\pi}\|\xi(\cdot;y^n_0)-\xi^*(\cdot)\|^2_{\mathcal{H}^2(\mathbb{C}^+_0;H)}
    \to 0\;\;\mbox{as}\;\;n\to+\infty.
\end{equation*}
    Thus,  for each $\lambda\in\mathbb{C}^+_0$,
\begin{eqnarray}\label{yu-5-12-6}
    \|\xi(\lambda;y^n_0)-\xi^*(\lambda)\|^2_H&\leq& \Big(\int_{\mathbb{R}^+}
    e^{-\mbox{Re}\lambda t}\|h^n(t)-h^*(t)\|_Hdt\Big)^2\nonumber\\
    &\leq&(2\mbox{Re}
    \lambda)^{-1}\|h^n-h^*\|^2_{L^2(\mathbb{R}^+;H)}\to0\;\;\mbox{as}\;\;n\to+\infty.
\end{eqnarray}
   Similarly, we can show that there is  $\eta^*(\cdot)\in \mathcal{H}^2(\mathbb{C}^+_0;U)$ such that
    for each $\lambda\in\mathbb{C}^+_0$,
\begin{equation}\label{yu-5-12-7}
    \|\eta(\lambda;y^n_0)-\eta^*(\lambda)\|_U\to 0\;\;\mbox{as}\;\;n\to+\infty.
\end{equation}
    Therefore, by (\ref{yu-5-12-1}), (\ref{yu-5-12-6}) and (\ref{yu-5-12-7}), we obtain that for $\varphi\in H_1$ and  $\lambda\in \mathbb{C}^+_0$,
\begin{eqnarray*}\label{yu-5-12-8}
    &\;&\langle \xi^*(\lambda),(\bar{\lambda}I-A^*)\varphi\rangle_H
    +\langle\eta^*(\lambda),B^*\varphi\rangle_U\nonumber\\
    &=&\lim_{n\to+\infty}\langle \xi(\lambda;y_0^n),(\bar{\lambda}I-A^*)\varphi\rangle_H
    +\langle\eta(\lambda;y_0^n),B^*\varphi\rangle_U
    =\lim_{n\to+\infty}\langle y_0^n,\varphi\rangle_H=\langle y_0,\varphi\rangle_H,
\end{eqnarray*}
which leads to (\ref{yu-5-12-1}). Meanwhile,  based on (\ref{yu-5-12-6}) and (\ref{yu-5-12-7}),
we can directly observe that  (\ref{yu-5-10-10}) holds for all $y_0\in H$.
Hence, Claim One has been proven.

\par
  Finally, we arbitrarily  fix  $y_0\in H$. Letting $\lambda=\beta+\mu$ with $\mu\in\mathbb{C}_{-\beta}^+$, formulas (\ref{yu-5-12-1}) and (\ref{yu-5-10-10}) lead to
\begin{equation}\label{yu-5-13-1}
    \langle \xi(\mu+\beta;y_0),(\bar{\mu}I-A^*)\varphi\rangle_H+\langle \eta(\mu+\beta;y_0),B^*\varphi\rangle_U=\langle y_0,\varphi\rangle_H
    \;\;\mbox{for}\;\;\mu\in\mathbb{C}_{-\beta}^+,\;\varphi\in H_1;
\end{equation}
    \begin{equation}\label{yu-5-13-2}
    \|\xi(\mu+\beta;y_0)\|_H\leq \frac{C_{\beta,\gamma,1}}{\mbox{Re}\,\mu+\beta}\|y_0\|_H
    \;\;\mbox{and}\;\;
    \|\eta(\mu+\beta;y_0)\|_U\leq \frac{C_{\beta,\gamma,2}}{\sqrt{\mbox{Re}\,\mu+\beta}}\|y_0\|_H\;\;
    \mbox{for}\;\;\mu\in\mathbb{C}^+_{-\beta}.
\end{equation}
   One can directly check  that $(\xi(\cdot+\beta;y_0),\eta(\cdot+\beta;y_0))
    \in \mathcal{H}^2(\mathbb{C}^+_{-\beta};H)\times \mathcal{H}^2(\mathbb{C}^+_{-\beta};U)$ (since $(\xi(\cdot;y_0),\eta(\cdot;y_0))
    \in \mathcal{H}^2(\mathbb{C}^+_{0};H)\times \mathcal{H}^2(\mathbb{C}^+_{0};U)$). Thus,  (\ref{yu-5-13-1}) and (\ref{yu-5-13-2}) imply that (\ref{yu-3-24-18}) and (\ref{yu-3-24-18-b}) hold with $C(\beta):=C_{\beta,\gamma,1}$ and $D(\beta):=C_{\beta,\gamma,2}$, respectively. Hence, conclusion
    $(ii)$ is true.

\vskip 5pt
  \noindent  \emph{Step 2. The proof of $(ii)\Rightarrow (iii)$ is trivial.}
\vskip 5pt
    \noindent \emph{Step 3. We prove  $(iii)\Rightarrow (i)$.}

    Suppose that $(iii)$ holds, i.e., there is $\beta\geq0$ such that for each $y_0\in H$, there  is $(\xi(\cdot;y_0),\eta(\cdot;y_0))\in \mathcal{H}^2(\mathbb{C}^+_{-\beta};H)\times \mathcal{H}^2(\mathbb{C}^+_{-\beta};U)$ satisfying (\ref{yu-3-24-18}). We arbitrarily fix $y_0\in H$. Then, based on $(i)$ in Lemma \ref{yu-lemma-6-9-1}, there is a unique $(y(\cdot),u(\cdot))\in L^2_{\beta}(\mathbb{R}^+;H)
    \times L^2_{\beta}(\mathbb{R}^+;U)$ such that
\begin{equation}\label{yu-5-14-1}
    \xi(\lambda;y_0)=\int_{\mathbb{R}^+}e^{-\lambda t}y(t)dt;\;\;\eta(\lambda;y_0)=
    -\int_{\mathbb{R}^+}e^{-\lambda t}u(t)dt,\; \lambda\in \mathbb{C}_{\beta}^+.
\end{equation}
    Let $\omega:=\max\left\{1,\lim_{t\to+\infty}t^{-1}\ln\|S(t)\|_{\mathcal{L}(H)}\right\}$.
    Then,    $\xi(\cdot;y_0)$ and $\eta(\cdot;y_0)$ are well defined over $\mathbb{C}_{\omega}^+$,
    and moreover, by \cite[Proposition 2.2, Chapter IV]{Engel-Nagel}, we have
   $\mathbb{C}^+_{\omega}\subset \rho(A)(=\rho(\widetilde{A}))$.
   We arbitrarily fix  $\lambda\in\mathbb{C}_{\omega}^+$. Then, we have that  $(\bar{\lambda} I-{A}^*)^{-1}\in\mathcal{L}(H;H_1)$, and  it follows from the proof of Theorem 3.1 in
   \cite[Section 1.3, Chapter 1]{Pazy} that
\begin{equation}\label{yu-5-14-2}
    (\bar{\lambda} I-A^*)^{-1}\varphi=\int_{\mathbb{R}^+}e^{-\bar{\lambda} t}S^*(t)\varphi\, dt
    \;\;\mbox{for each}\;\;\varphi\in H.
\end{equation}
       We arbitrarily fix  $\psi\in H_1$.
      There are two facts. First, replacing $\varphi$ by $(\bar{\lambda}I-A^*)^{-1}\psi$ in  (\ref{yu-3-24-18}) ($\lambda\in\mathbb{C}_{\omega}^+$) leads to
\begin{equation}\label{yu-5-14-3}
    \langle\xi(\lambda;y_0),\psi\rangle_H+
    \langle\eta(\lambda;y_0),B^*(\bar{\lambda}I-A^*)^{-1}\psi\rangle_U=\langle y_0,(\bar{\lambda}I-A^*)^{-1}\psi\rangle_H.
\end{equation}
   Second, with assumption $(A_2)$, we have
   $B^*\in \mathcal{L}(H_1;U)$ and
   $B^*(\lambda I-A^*)^{-1}\in \mathcal{L}(H;U)$, thus,
\begin{equation}\label{yu-5-14-4}
    \langle\eta(\lambda;y_0),B^*(\bar{\lambda}I-A^*)^{-1}\psi\rangle_U
    =\langle B\eta(\lambda;y_0),(\bar{\lambda}I-A^*)^{-1}\psi\rangle_{H_{-1},H_1}.
   \end{equation}

\par

   Now, we claim that
\begin{equation}\label{yu-5-14-b-1}
    e^{-\mbox{Re}\lambda \cdot}\int_{0}^{\cdot} \widetilde{S}(\cdot-s)Bu(s)ds\in L^1(\mathbb{R}^+;H_{-1})\cap L^2(\mathbb{R}^+;H_{-1}).
\end{equation}
    For this purpose, we first recall that $\rho_0$ is provided in assumption $(A_2)$, so we have $(\rho_0I-\widetilde{A})^{-1}B\in \mathcal{L}(U;H)$. Moreover,
    by \cite[Proposition 2.10.3, Chapter 2]{Tucsnak-Weiss}, we have that
     $\rho_0I-\widetilde{A}$ is a unitary linear map from $H$ to $H_{-1}$. Hence,   for each $t\in \mathbb{R}^+$,
\begin{eqnarray*}
    &\;&e^{-\mbox{Re}\lambda t}\Big\|\int_0^t\widetilde{S}(t-s)Bu(s)ds\Big\|_{H_{-1}}
    =e^{-\mbox{Re}\lambda t}\Big\|\int_0^tS(t-s)(\rho_0I-\widetilde{A})^{-1}Bu(s)ds\Big\|_H\nonumber\\
    &\leq& C(\omega)\|(\rho_0I-\widetilde{A})^{-1}B\|_{\mathcal{L}(U;H)}
    e^{-\mbox{Re}\lambda t}\int_0^t e^{\omega (t-s)}\|u(s)\|_Uds\nonumber\\
    &\leq&C(\omega)\|(\rho_0I-\widetilde{A})^{-1}B\|_{\mathcal{L}(U;H)}
    e^{-\mbox{Re}\lambda t}\Big(\int_0^te^{2\omega s}ds\Big)^{\frac{1}{2}}
    \Big(\int_0^t\|u(s)\|_U^2ds\Big)^{\frac{1}{2}}\nonumber\\
    &\leq&C(\omega)\|(\rho_0I-\widetilde{A})^{-1}B\|_{\mathcal{L}(U;H)}(2\omega)^{-\frac{1}{2}}
    \|u\|_{L^2_\beta(\mathbb{R}^+;U)} e^{-(\mbox{Re}\lambda-\omega)t}.
\end{eqnarray*}
Since $\mbox{Re}\lambda>\omega$, the above leads to  (\ref{yu-5-14-b-1}).

Based on $(i)$ in Remark \ref{yu-b-remark-5-11-1}, (\ref{yu-5-14-1}), (\ref{yu-5-14-2})
and (\ref{yu-5-14-b-1}),
using the Fubini theorem and $u=0$ in $(-\infty,0)$, we obtain
\begin{eqnarray*}\label{yu-5-14-5}
     &\;&\langle B\eta(\lambda;y_0),(\bar{\lambda}I-A^*)^{-1}\psi\rangle_{H_{-1},H_1}\nonumber\\
     &=&
     -\Big\langle B\int_{\mathbb{R}^+}e^{-\lambda t}u(t)dt,\int_{\mathbb{R}^+}
     e^{-\bar{\lambda}t}S^*(t)\psi dt\Big\rangle_{H_{-1},H_1}
     \nonumber\\
     &=&-\Big\langle (\rho_0 I-\widetilde{A})^{-1}B\int_{\mathbb{R}^+}e^{-\lambda t}u(t)dt,(\rho_0 I-A^*)\int_{\mathbb{R}^+}
     e^{-\bar{\lambda}t}S^*(t)\psi dt\Big\rangle_H\nonumber\\
     &=&-\int_{\mathbb{R}^+}\int_{t}^{+\infty}\langle e^{-\lambda(\sigma-t)}(\rho_0 I-\widetilde{A})^{-1}Bu(\sigma-t), e^{-\bar{\lambda}t}S^*(t)(\rho_0 I-A^*)\psi\rangle_{H}d\sigma dt\nonumber\\
     &=&-\int_{\mathbb{R}^+}e^{-\lambda \sigma}\int_{0}^{\sigma}\langle
     (\rho_0 I-\widetilde{A})^{-1}Bu(\sigma-t),S^*(t)(\rho_0 I-A^*)\psi\rangle_{H}dt d\sigma\nonumber\\
     &=&-\int_{\mathbb{R}^+}e^{-\lambda \sigma}\Big\langle \int_{0}^\sigma \widetilde{S}(\sigma-t)Bu(t)dt,\psi
     \Big\rangle_{H{-1},H_1}d\sigma.
\end{eqnarray*}
     The combination of this result, (\ref{yu-5-14-3}) and (\ref{yu-5-14-4}) imply that
    \begin{equation}\label{w2.22-7-20}
     \int_{\mathbb{R}}e^{-\lambda t}
    \left\langle F(t),\psi\right\rangle_{H_{-1},H_1}dt=0\;\;\mbox{for} \;\;\lambda\in \mathbb{C}_{\omega}^+,\;\; \psi\in H_1,
     \end{equation}
     where $F(t):=y(t)-\widetilde{S}(t)y_0-\int_0^t\widetilde{S}(t-s)Bu(s)ds$ for $t\geq 0$, while
     $F(t):=0$ for $t<0$.
                        With (\ref{yu-5-14-b-1}),  we can apply the  inverse Fourier transform to (\ref{w2.22-7-20}) with respect to $\mbox{Im}\lambda$ to  conclude that
\begin{equation*}\label{yu-5-14-7}
    y(t)=\widetilde{S}(t)y_0+\int_0^t\widetilde{S}(t-s)Bu(s)ds\;\;\mbox{a.e.}\;\;t\in\mathbb{R}^+,
\end{equation*}
which leads to
   $y(\cdot)=y(\cdot;y_0,u)$ a.e. in $\mathbb{R}^+$.
    With
     $(y(\cdot),u(\cdot))\in L^2_{\beta}(\mathbb{R}^+;H)
    \times L^2_{\beta}(\mathbb{R}^+;U)\subset L^2(\mathbb{R}^+;H)
    \times L^2(\mathbb{R}^+;U)$, we obtain
    \begin{equation}\label{w2.24-7-20}
    \mathcal{U}_{ad}(y_0):=\{u(\cdot)\in L^2(\mathbb{R}^+;U):y(\cdot;y_0,u)\in L^2(\mathbb{R}^+;H)\}\neq\emptyset.
    \end{equation}
        Since $y_0$ was arbitrarily taken from $H$, \eqref{w2.24-7-20}
        and
     \cite[Proposition 3.9]{Liu-Wang-Xu-Yu} imply that system $[A,B]$ is exponentially stabilizable.
   Hence, we have completed the proof of Theorem \ref{yu-main-theorem-1}.
\end{proof}

\section{Proofs of the main theorems}
Before presenting the proofs of the main theorems,  we need the following lemmas:
\begin{lemma}\label{yu-lemma-7-6-1}
    Suppose that $A\in \mathcal{L}(H)$ and $B\in\mathcal{L}(U;H)$. Then, the following statements are equivalent:
\begin{enumerate}
  \item [$(i)$] System $[A,B]$ is exactly controllable at some $T>0$;
  \item[$(ii)$] There are constants $T>0$ and $C(T)>0$ such that
\begin{equation*}
    \|\varphi\|_H^2\leq C(T)\int_0^T\|B^*e^{A^*t}\varphi\|_U^2dt\;\;\mbox{for}\;\;\varphi\in H;
\end{equation*}
  \item [$(iii)$] There is $n\in\mathbb{N}$ such that
  $\mbox{Span}\{BU,ABU,\ldots,A^nBU\}=H$;
  \item [$(iv)$] For each $\lambda\in\mathbb{C}$, there is $C(\lambda)>0$ such that $\|\varphi\|_H^2
  \leq C(\lambda)(\|(\lambda I-A^*)\varphi\|_H^2+\|B^*\varphi\|^2_U)$ for all $\varphi\in H$.
\end{enumerate}
\end{lemma}
\begin{proof}
    It is well-known that  $(i)\Leftrightarrow (ii)$ (e.g.,
    \cite[Theorem 11.2.1, Chapter 11]{Tucsnak-Weiss}).
    Using Baire category theorem, one can directly verify that $(i)\Rightarrow (iii)$, while
    the proof of $(iii)\Rightarrow (i)$ can be found in \cite[Theorem 2.3]{Triggiani}.
    The proof of $(iii)\Leftrightarrow (iv)$ can be found in \cite{Takahashi}
    (see the main theorem and the remark on it there). This completes the proof.
\end{proof}

\begin{lemma}\label{lemma3.2-7-22}
System $[A,B]$ is rapidly stabilizable if and only if for each $\alpha>0$, system $[A_\alpha,B]$ is  stabilizable, where $A_\alpha:=A+\alpha I$.
\end{lemma}
 If $B\in\mathcal{L}(U;H)$, then Lemma \ref{lemma3.2-7-22}  is well known. However,  for our framework,  we do not find an accurate literature with its proof. Thus, we prove it here.
\begin{proof}[Proof of Lemma \ref{lemma3.2-7-22}.]
First, let  $S_\alpha(\cdot)$ be the $C_0$-semigroup generated by $A_\alpha$. Then, we have $S_\alpha(t)=e^{\alpha t}S(t)$, $t\geq 0$.

Now, we suppose that  system $[A,B]$ is rapidly stabilizable. Then, according to
    \cite[Theorem 3.4]{Liu-Wang-Xu-Yu},  for each $\alpha>0$, there is $C(\alpha)>0$ such that
\begin{equation*}\label{yu-7-2-10}
    \|S^*(t)\varphi\|_H^2\leq C(\alpha)\left(\int_0^t\|B^*S^*(s)\varphi\|^2_Uds
    +e^{-(2\alpha+1)t}\|\varphi\|_H^2\right)\;\;\mbox{for}\;\;\varphi\in H_1,\;t\in\mathbb{R}^+.
\end{equation*}
    Hence, for $\alpha>0$,
\begin{equation}\label{yu-7-2-11}
    \|S^*_\alpha(t)\varphi\|_H^2\leq C(\alpha)e^{2\alpha t}\int_0^t\|B^*S^*_\alpha(s)\varphi\|_U^2ds
    +C(\alpha)e^{-t}\|\varphi\|_H^2\;\;\mbox{for}\;\;\varphi\in H_1,\;t\in\mathbb{R}^+.
\end{equation}
  For each  $\alpha>0$, we let  $T>0$ satisfy $\delta:=D(\alpha)e^{-T}<1$. Then, based on (\ref{yu-7-2-11}), we have
\begin{equation*}\label{yu-7-2-12}
    \|S^*_\alpha(T)\varphi\|_H^2\leq C_1(\alpha)e^{2\alpha T}\int_0^T
    \|B^*S^*_{\alpha}(s)\varphi\|_U^2ds+\delta\|\varphi\|_H^2\;\;\mbox{for}\;\;
    \varphi\in H_1.
    \end{equation*}
    The combination of this and $(i)$ of Lemma \ref{yu-bb-remark-6-21-1} yield that
    system $[A_\alpha, B]$ is  stabilizable.
\par
    Conversely, we suppose that,  for each $\alpha>0$, system $[A_\alpha,B]$ is  stabilizable. Then, according to   $(i)$ of Lemma \ref{yu-bb-remark-6-21-1}, for each $\alpha>0$, there is $C(\alpha)>0$ such that
$$
    \|S^*_\alpha(t)\varphi\|_H^2\leq C(\alpha)\left(\int_0^t\|B^*S^*_\alpha(s)\varphi\|_U^2ds
    +\|\varphi\|_H^2\right)\;\;\mbox{for}\;\;\varphi\in H_1,\;t\in\mathbb{R}^+,
$$
  which implies that  if $\alpha>0$,
$$
    \|S^*(t)\varphi\|_H^2\leq C(\alpha)\left(\int_0^t\|B^*S^*(s)\varphi\|_U^2ds
    +e^{-\alpha t}\|\varphi\|_H^2\right)\;\;\mbox{for}\;\;\varphi\in H_1,\;t\in\mathbb{R}^+.
$$
     With $(ii)$ of  Lemma \ref{yu-bb-remark-6-21-1}, the above leads to the rapid stabilizability of system
     $[A,B]$.
     This completes the proof of Lemma \ref{lemma3.2-7-22}.
\end{proof}
    We start with proving Theorem \ref{yu-theorem-7-3-1}.
\begin{proof}[Proof of Theorem \ref{yu-theorem-7-3-1}.]
 We first prove $(i)$. Suppose that system $[A,B]$ is  stabilizable.
    According to Theorem \ref{yu-main-theorem-1},  there are constants $\beta_0>0$, $C(\beta_0)\geq 1$ and $D(\beta_0)>0$ such that for each $y_0\in H$, there is $(\xi(\cdot;y_0),\eta(\cdot;y_0))
    \in \mathcal{H}^2(\mathbb{C}_{-\beta_0}^+;H)\times \mathcal{H}^2(\mathbb{C}_{-\beta_0}^+;U)$ that satisfies (\ref{yu-3-24-18}) and (\ref{yu-3-24-18-b}). Then, based on (\ref{yu-3-24-18}) and (\ref{yu-3-24-18-b}), for $\lambda\in\mathbb{C}_{-\beta_0}^+$ and $y_0\in H$,
 \begin{eqnarray*}\label{yu-5-15-20}
    |\langle y_0,\varphi\rangle_H|\leq
    \Big(\frac{C(\beta_0)}{\mbox{Re}\,\lambda+\beta_0}\|(\bar{\lambda} I-A^*)\varphi\|_H
    +\frac{D(\beta_0)}{\sqrt{\mbox{Re}\,\lambda+\beta_0}}\|B^*\varphi\|_U\Big)\|y_0\|_H
    \;\;\mbox{for}\;\;\varphi\in H_1.
\end{eqnarray*}
    Thus, for each $\lambda\in\mathbb{C}_{-\beta_0}^+$,
\begin{equation}\label{yu-5-16-1}
    \|\varphi\|_{H}\leq \frac{C(\beta_0)}{\mbox{Re}\,\lambda+\beta_0}\|(\bar{\lambda} I-A^*)\varphi\|_H
    +\frac{D(\beta_0)}{\sqrt{\mbox{Re}\,\lambda+\beta_0}}\|B^*\varphi\|_U
    \;\;\mbox{for}\;\;\varphi\in H_1.
\end{equation}
   Suppose that
        $\|S^*(t)\|\leq C(\omega)e^{\omega t}$ for each $t\in\mathbb{R}^+$ for some constants $\omega>0$ and $C(\omega)>0$. Taking $\beta\in(0,\beta_0)$.
        By \cite[Theorem 5.3 and Remark 5.4, Section 1.5, Chapter 1] {Pazy}, we have that, for each
        $\lambda\in \mathbb{C}^+_{\max\{\omega,2|\beta-\omega|-\beta\}}$,
\begin{eqnarray}\label{yu-8-17-1}
    \|\varphi\|_H&\leq& \frac{C(\omega)}{\mbox{Re}\,\lambda-\omega}
    \|(\lambda I-A^*)\varphi\|_H\leq \frac{C(\omega)}{(\mbox{Re}\,\lambda+\beta)
    -|\beta-\omega|}
    \|(\lambda I-A^*)\varphi\|_H\nonumber\\
    &=&\frac{C(\omega)}{\frac{1}{2}(\mbox{Re}\,\lambda+\beta)+(\frac{1}{2}(\mbox{Re}\,\lambda
    +\beta)
    -|\beta-\omega|)}
    \|(\lambda I-A^*)\varphi\|_H\nonumber\\
    &\leq& \frac{2C(\omega)}{\mbox{Re}\,\lambda+\beta}
    \|(\lambda I-A^*)\varphi\|_H \;\;\mbox{for}
    \;\;\varphi\in H_1.
\end{eqnarray}
       If $\lambda\in \mathbb{C}_{-\beta}^+\setminus \mathbb{C}^+_{\max\{\omega,2|\beta-\omega|-\beta\}}$, then
\begin{equation*}
    \frac{D(\beta)}{\sqrt{\mbox{Re}\,\lambda+\beta_0}}
    \leq \frac{D(\beta)(\max\{\omega,2|\beta-\omega|-\beta\}+\beta)}{\sqrt{\beta_0-\beta}}
    \frac{1}{\mbox{Re}\,\lambda+\beta}.
\end{equation*}
   The combination of this relation, \eqref{yu-8-17-1}, \eqref{yu-5-16-1} and Definition \ref{yu-definition-6-9-1} imply that
      the pair $(A^*,B^*)$ satisfies (HESI)$_\beta$.
\par
    Next, we prove $(ii)$. Suppose that system $[A,B]$ is rapidly stabilizable. Lemma \ref{lemma3.2-7-22} implies that, for an arbitrarily fixed  $\alpha>0$, system $[A+\alpha I,B]$ is stabilizable. Then, according to the conclusion  $(i)$ of Theorem \ref{yu-theorem-7-3-1},  there is $C(\alpha)>0$ such that if $\lambda\in\mathbb{C}_0^+$,
$$
    \|\varphi\|_H^2\leq \frac{C(\alpha)}{(\mbox{Re}\,\lambda)^2}\left(\|(\lambda I-(A^*+\alpha I))\varphi
    \|_H^2+\|B^*\varphi\|_U^2\right)\;\;\mbox{for}\;\;\varphi\in H_1.
$$
    Hence, if $\lambda\in\mathbb{C}^+_{-\alpha}$,
$$
    \|\varphi\|_H^2\leq \frac{C(\alpha)}{(\mbox{Re}\,\lambda+\alpha)^2}\left(\|(\lambda I-A^*)\varphi
    \|_H^2+\|B^*\varphi\|_U^2\right)\;\;\mbox{for}\;\;\varphi\in H_1.
$$
    Since  $\alpha>0$ was arbitrarily taken, the above shows that  the pair $(A^*,B^*)$ satisfies (HESI).
       This completes the proof of Theorem \ref{yu-theorem-7-3-1}.
              \end{proof}

We now prove Theorem \ref{yu-theorem-6-30-1}.
\begin{proof}[Proof of Theorem \ref{yu-theorem-6-30-1}.]
First, since $A^*$ is normal, we can write $E^{A^*}$ for the unique spectral measure corresponding to $A^*$, which is provided by
the   spectral theorem (e.g., \cite[Theorem 6.47, Chapter 6]{Kubrusly}).
    We divide the proof into two steps.
\vskip 5pt
   \noindent \emph{Step 1. We prove conclusion $(i)$.}

   By Theorem
    \ref{yu-theorem-7-3-1}, we have the necessity.
    The remainder is to show the sufficiency.
   We suppose that
   the pair $(A^*,B^*)$ satisfies (HESI)$_\beta$ for some $\beta>0$. Then, there is $\beta>0$ and $C(\beta)>0$ satisfying  (\ref{yu-6-9-1}). Without loss of generality, we can assume that  $\mathbb{C}_{-\varepsilon}^+\cap \sigma(A^*)\neq \emptyset$ for each
    $\varepsilon>0$. (Otherwise,  there is $\varepsilon^*>0$ such that $\mathbb{C}_{-\varepsilon^*}^+\cap \sigma(A^*)=\emptyset$. Then, it follows from
    \cite[Corollary 3.4, Section 3, Chapter IV and Lemma 1.9, Section 1, Chapter V]{Engel-Nagel} that
    $S^*(\cdot)$ and $S(\cdot)$ are exponentially stable. Thus, by taking the feedback law as $0$, we obtain the sufficiency.)
     We take $\omega>0$ such that $\|S^*(t)\|_{\mathcal{L}(H)}\leq C(\omega)e^{\omega t}$ for each $t\in\mathbb{R}^+$. The remainder of the  proof in this step is organized into two sub-steps.
\vskip 5pt
 \noindent   \emph{Sub-step 1.1. We prove that for each $\beta^*\in(0,\beta)$, there are $T:=T(\beta^*)>0$ and $C(T,\beta^*)>0$ such that
\begin{equation}\label{yu-7-12-1}
    \|S^*(T)E^{A^*}({\mathbb{C}_{-{\beta^*}}^+})
    \varphi\|_H^2\leq C(T,\beta^*)
    \int_0^{T}\|B^*E^{A^*}({\mathbb{C}_{-{\beta^*}}^+})
    S^*(t)E^{A^*}({\mathbb{C}_{-{\beta^*}}^+})\varphi\|_U^2dt
    \;\;\mbox{for}\;\;\varphi\in H_1.
\end{equation} }
    We arbitrarily fix $\beta^*\in(0,\beta)$. Since $\sigma(A^*)\cap \mathbb{C}_{-\gamma}^+$ is bounded
     (by our assumption),  it follows from  the spectral theorem that
\begin{enumerate}
  \item [$(a)$] $A^*E^{A^*}({\mathbb{C}_{-\beta^*}^+})=E^{A^*}({\mathbb{C}_{-\beta^*}^+})A^*$ is a bounded operator on $H$;
  \item [$(b)$] $E^{A^*}({\mathbb{C}_{-\beta^*}^+})H$ is an invariant subspace of $A^*E^{A^*}({\mathbb{C}_{-\beta^*}^+})$.
\end{enumerate}
    Thus, we have  $E^{A^*}({\mathbb{C}_{-\beta^*}^+})H= D(A^*E^{A^*}({\mathbb{C}_{-\beta^*}^+}))$
    and $A^*E^{A^*}({\mathbb{C}_{-\beta^*}^+})\in\mathcal{L}(E^{A^*}({\mathbb{C}_{-\beta^*}^+})H)$.
    Consequently,
     $E^{A^*}({\mathbb{C}_{-\beta^*}^+})S^*(\cdot)(=S^*(\cdot)E^{A^*}({\mathbb{C}_{-\beta^*}^+}))$ is the $C_0$-semigroup on $H$, which is generated by  $A^*E^{A^*}({\mathbb{C}_{-\beta^*}^+})$.
    Based on $(b)$, we know that  $E^{A^*}({\mathbb{C}_{-\beta^*}^+})H$ is an invariant subspace of $S^*(t)$ for each $t\in\mathbb{R}^+$, so $E^{A^*}({\mathbb{C}_{-\beta^*}^+})S^*(\cdot)$ is  a $C_0$-semigroup on $E^{A^*}({\mathbb{C}_{-\beta^*}^+})H$. Meanwhile, based on $(a)$ and assumption $(A_2)$ (see also
     $(i)$ in Remark \ref{yu-b-remark-5-11-1}),
     $B^*E^{A^*}({\mathbb{C}_{-\beta^*}^+})\in \mathcal{L}(E^{A^*}({\mathbb{C}_{-\beta^*}^+})H;U)$.
     Indeed, we have
$$
    \|B^*E^{A^*}({\mathbb{C}_{-\beta^*}^+})\varphi\|_U
    \leq \|(B^*(\rho_0I-A^*)^{-1})\|_{\mathcal{L}(H;U)}
    \|(\rho_0I-A^*)E^{A^*}({\mathbb{C}_{-\beta^*}^+})\|_{\mathcal{L}(H)}
    \|\varphi\|_H
$$
    for all $\varphi\in E^{A^*}({\mathbb{C}_{-\beta^*}^+})H$.
\par
    Now, we claim that there is $C(\beta,\beta^*)>0$ such that if $\lambda\in \mathbb{C}$,
\begin{equation}\label{yu-7-6-10}
    \|\varphi\|_H^2\leq C(\beta,\beta^*)(\|(\lambda I-A^*E^{A^*}({\mathbb{C}_{-\beta^*}^+}))\varphi\|^2_H
    +\|B^*E^{A^*}({\mathbb{C}_{-\beta^*}^+})\varphi\|_U^2)\;\;\mbox{for}
    \;\;\varphi\in E^{A^*}({\mathbb{C}_{-\beta^*}^+})H.
\end{equation}
     Indeed, (\ref{yu-6-9-1}) implies that for $\lambda\in \overline{\mathbb{C}^+_{-\frac{\beta^*+\beta}{2}}}
     (\subset \mathbb{C}_{-\beta}^+)$ and $\varphi\in E^{A^*}({\mathbb{C}_{-\beta^*}^+})H$,
\begin{equation}\label{yu-7-6-11}
    \|\varphi\|_H^2\leq
    \frac{4C(\beta)}{(\beta-\beta^*)^{2}}
    (\|(\lambda I-A^*E^{A^*}(\mathbb{C}_{-\beta^*}^+))\varphi\|^2_H
    +\|B^*E^{A^*}({\mathbb{C}_{-\beta^*}^+})\varphi\|_U^2).
\end{equation}
     However,
     according to  the spectral theorem,
\begin{equation*}\label{yu-7-6-12}
    \|(\lambda I-A^*E^{A^*}({\mathbb{C}_{-\beta^*}^+}))^{-1}\varphi
    \|_H^2\leq \frac{4}{(\beta-\beta^*)^2}\|\varphi\|_H^2\;\; \mbox{for}\;\;\lambda\in \mathbb{C}^-_{-\frac{\beta^*+\beta}{2}},\;  \varphi\in E^{A^*}({\mathbb{C}_{-\beta^*}^+})H,
\end{equation*}
   which implies that
\begin{equation*}\label{yu-7-6-13}
    \|\varphi
    \|_H^2\leq \frac{4}{(\beta-\beta^*)^2}\|(\lambda I-A^*E^{A^*}({\mathbb{C}_{-\beta^*}^+}))\varphi\|_H^2
    \;\;\mbox{for}\;\;\lambda\in \mathbb{C}^-_{-\frac{\beta^*+\beta}{2}},\;\varphi\in E^{A^*}({\mathbb{C}_{-\beta^*}^+})H.
\end{equation*}
     This fact and (\ref{yu-7-6-11}) lead to (\ref{yu-7-6-10})  with
     $C(\beta,\beta^*):=\frac{4(1+C(\beta))}{(\beta-\beta^*)^2}$.
     Since
     $$
     (A^*E^{A^*}({\mathbb{C}_{-\beta^*}^+}),B^*E^{A^*}({\mathbb{C}_{-\beta^*}^+}))
     \in\mathcal{L}(E^{A^*}({\mathbb{C}_{-\beta^*}^+})H)
     \times \mathcal{L}(E^{A^*}({\mathbb{C}_{-\beta^*}^+})H;U),
     $$
     (\ref{yu-7-6-10}) and Lemma \ref{yu-lemma-7-6-1} imply that
      there is $T:=T(\beta^*)>0$ and $C(T,\beta^*)>0$ that satisfies  (\ref{yu-7-12-1}).
     (Here, we use that $\|S^*(T)E^{A^*}({\mathbb{C}_{-{\beta^*}}^+})
    \varphi\|_H\leq C(\omega)e^{\omega T}\|E^{A^*}({\mathbb{C}_{-{\beta^*}}^+})
    \varphi\|_H$
    for each $\varphi\in H_1$.) Thus, we have completed {\it Sub-step 1.1}.

\vskip 5pt
   \noindent \emph{Sub-step 1.2. We prove that the system $[A,B]$ is  stabilizable.}

   We arbitrarily fix $\beta^*\in(0,\beta)$. The facts are as follows:
   First, according to \emph{Substep 1.1}, there are constants
      $T_0>0$ and  $C(T_0, \beta^*)>0$ such that (\ref{yu-7-12-1}) (where $T=T_0$ and $C(T,\beta^*)=C(T_0,\beta^*)$) holds.
      Second, one can easily check  that $A^*E^{A^*}(\overline{\mathbb{C}_{-{\beta^*}}^-})$ (with its domain
       $E^{A^*}(\overline{\mathbb{C}_{-{\beta^*}}^-})H_1$) generates the $C_0$-semigroup
            $S^*(\cdot) E^{A^*}(\overline{\mathbb{C}_{-{\beta^*}}^-})$  on $E^{A^*}(\overline{\mathbb{C}_{-{\beta^*}}^-})H$. Third, the spectral theorem implies that $\sigma(A^*E^{A^*}(\overline{\mathbb{C}_{-{\beta^*}}^-}))\subset \overline{\mathbb{C}_{-{\beta^*}}^-}$ and thus
\begin{equation}\label{yu-7-2-2}
    \sup\{\mbox{Re}\lambda:\lambda\in \sigma(A^*E^{A^*}(\overline{\mathbb{C}_{-{\beta^*}}^-}))\}\leq -\beta^*.
\end{equation}
    Since $A^*$ is a normal operator and $E^{A^*}(\overline{\mathbb{C}_{-{\beta^*}}^-})$ is an orthogonal projection, we have that $A^*E^{A^*}(\overline{\mathbb{C}_{-{\beta^*}}^-})$ is  normal on $E^{A^*}(\overline{\mathbb{C}_{-{\beta^*}}^-})H$. Then, by (\ref{yu-7-2-2}) and
    \cite[Corollary 3.4, Section 3, Chapter IV and Lemma 1.9, Section 1, Chapter V]{Engel-Nagel}, we obtain that for each
    $\eta\in(0,\beta^*)$, there is $C(\eta)>0$ such that
\begin{equation}\label{yu-7-2-3}
    \|S^*(t)E^{A^*}(\overline{\mathbb{C}_{-{\beta^*}}^-})\|_{\mathcal{L}(H)}\leq
    C(\eta)e^{-\eta t}\;\;\mbox{for}\;\;t\in\mathbb{R}^+.
\end{equation}
\par
    Now, we claim that for each $T\geq 2T_0$, there is $C(T_0, \omega,\beta^*,\eta)>0$ such that
\begin{equation}\label{yu-7-2-4}
    \|S^*(T)\varphi\|_H^2\leq C(T_0, \omega,\beta^*,\eta)
    \Big(\int_0^T\|B^*S^*(t)\varphi\|_U^2dt+e^{-2\eta T}\|\varphi\|_H^2\Big)
    \;\;\mbox{for}\;\;\varphi\in H_1.
\end{equation}
    For this purpose, we arbitrarily fix $T\geq 2T_0$. Let $N:=[T/T_0]$,
    then $N\geq 2$ and $NT_0\leq T<(N+1)T_0$. Based on (\ref{yu-7-12-1}) (with $T=T_0$ and $C(T,\beta^*)=C(T_0,\beta^*)$) and (\ref{yu-7-2-3}), we have that,
    for each $\varphi\in H_1$,
\begin{eqnarray}\label{yu-7-2-5}
    &&\|S^*(T)\varphi\|_H^2=\|S^*(T-NT_0)S^*(NT_0)\varphi\|_H^2
    \leq (C(\omega))^2e^{2\omega T_0}\|S^*(NT_0)\varphi\|_H^2\nonumber\\
    &=& (C(\omega))^2e^{2\omega T_0}
    \Big(\|S^*(T_0)E^{A^*}({\mathbb{C}_{-\beta^*}^+})S^*((N-1)T_0)\varphi\|^2_H
    +\|S^*(NT_0)E^{A^*}(\overline{\mathbb{C}_{-{\beta^*}}^-})\varphi\|_H^2\Big)\nonumber\\
    &\leq&(C(\omega))^2e^{2\omega T_0}
    \Big(C(T_0,\beta^*)\int_0^{T_0}\|B^*S^*((N-1)T_0+t)E^{A^*}
    ({\mathbb{C}_{-\beta^*}^+})\varphi\|_U^2dt
    +\|S^*(NT_0)E^{A^*}(\overline{\mathbb{C}_{-{\beta^*}}^-})\varphi\|_H^2\Big)\nonumber\\
    &\leq&I_1+I_2,
\end{eqnarray}
    where
\begin{equation*}\label{yu-7-2-6}
    I_1:=(C(\omega))^2e^{2\omega T_0}\Big(2C(T_0,\beta^*)\int_{(N-1)T_0}^{NT_0}
    \|B^*S^*(t)\varphi\|_U^2dt
    +(C(\eta))^2e^{-2\eta N T_0}\|\varphi\|_H^2\Big),
\end{equation*}
    and
\begin{equation*}\label{yu-7-2-7}
    I_2:=2(C(\omega))^2C(T_0,\beta^*)e^{2\omega T_0}
    \int_{0}^{T_0}\|B^*S^*((N-1)T_0+t)E^{A^*}(\overline{\mathbb{C}_{-{\beta^*}}^-})\varphi\|_H^2dt.
\end{equation*}
    Based on assumption $(A_3)$ and (\ref{yu-7-2-3}), one can directly check  that
\begin{eqnarray}\label{yu-7-2-8-b}
    I_1\leq 2(C(\omega))^2 C(T_0,\beta^*)e^{2\omega T_0}\int_0^T\|B^*S^*(t)\varphi\|_U^2dt
    +(C(\omega)C(\eta))^2e^{2(\omega+\eta)T_0}e^{-2\eta T}\|\varphi\|^2_H;
\end{eqnarray}
    \begin{eqnarray*}\label{yu-7-2-8}
    I_2
    &=&2(C(\omega))^2C(T_0,\beta^*)e^{2\omega T_0}\int_{0}^{T_0}\|B^*S^*(t)S^*((N-1)T_0)E^{A^*}(\overline{\mathbb{C}_{-{\beta^*}}^-})
    \varphi\|_H^2dt\nonumber\\
    &\leq& 2(C(\omega))^2C(T_0,\beta^*)C(T_0)e^{2\omega T_0}\|S^*((N-1)T_0)E^{A^*}(\overline{\mathbb{C}_{-{\beta^*}}^-})\varphi\|_H^2\nonumber\\
    &\leq& 2(C(\omega))^2C(T_0,\beta^*)C(T_0)(C(\eta))^2 e^{2\omega T_0}e^{-2\eta (N-1)T_0}\|E^{A^*}(\overline{\mathbb{C}_{-{\beta^*}}^-})\varphi\|_H^2\nonumber\\
    &\leq& 2(C(\omega)C(\eta))^2C(T_0,\beta^*)C(T_0)e^{2(\omega+2\eta)T_0}e^{-2\eta T}\|\varphi\|_H^2.
\end{eqnarray*}
     The above, (\ref{yu-7-2-5}), and (\ref{yu-7-2-8-b}) lead to (\ref{yu-7-2-4})  with
$$
    C(T_0, \omega,\beta^*,\eta):=(C(\omega))^2e^{2\omega T_0}
    \max\{2 C(T_0,\beta^*), (C(\eta))^2e^{2\eta T_0}(1+2C(T_0,\beta^*)C(T_0)e^{2\eta T_0})\}.
$$

     Using (\ref{yu-7-2-4}), we can find $\widehat{T}>0$ such that
\begin{equation*}\label{yu-7-2-9}
    \|S^*(\widehat{T})\varphi\|_H^2\leq C(T_0, \omega,\beta^*,\eta)
    \int_0^{\widehat{T}}\|B^*S^*(t)\varphi\|_U^2dt+\frac{1}{2}\|\varphi\|_H^2
    \;\;\mbox{for}\;\;\varphi\in H_1.
\end{equation*}
    The combination of this result and Lemma \ref{yu-bb-remark-6-21-1} yield that system $[A,B]$ is  stabilizable. Thus, we have completed {\it Sub-step 1.2} and {\it Step 1}.
\vskip 5pt
   \noindent \emph{Step 2. We prove conclusion $(ii)$.}
\par
    The necessity is proven in Theorem \ref{yu-theorem-7-3-1}. Thus, we must only prove the sufficiency.
   Suppose that  the pair $(A^*,B^*)$ satisfies (HESI). Then, for each $\beta>0$, there is $C(\beta)>0$ such that
\begin{equation*}\label{yu-7-2-16}
    \|\varphi\|^2_H\leq \frac{C(\beta)}{(\mbox{Re}\,\lambda+2\beta)^2}
    \left(\|(\lambda I-A^*)\varphi\|_H^2+\|B\varphi\|_U^2\right)
    \;\;\mbox{for}\;\;\lambda\in\mathbb{C}_{-2\beta}^+,\;\varphi\in H_1,
\end{equation*}
    which yields
\begin{equation*}\label{yu-7-2-17}
    \|\varphi\|^2_H\leq \frac{C(\beta)}{(\mbox{Re}\,\lambda+\beta)^2}
    \left(\|(\lambda I-(A^*+\beta I))\varphi\|_H^2+\|B\varphi\|_U^2\right)
    \;\;\mbox{for}\;\;\lambda\in\mathbb{C}_{-\beta}^+,\;\varphi\in H_1.
\end{equation*}
    The combination of this result and conclusion $(i)$ of Theorem \ref{yu-theorem-6-30-1} imply that
    $[A+\beta I,B]$ is  stabilizable. (Since $A$ is normal, $A+\beta I$ is also normal for any $\beta\in\mathbb{R}$.)  Since  $\beta>0$ can be arbitrarily chosen, Lemma \ref{lemma3.2-7-22} implies that $[A,B]$ is rapidly stabilizable.
\par
Thus, we complete the proof of Theorem \ref{yu-theorem-6-30-1}.
\end{proof}

Finally, we prove Theorem \ref{yu-theorem-5-20-1}.
\begin{proof}[The proof of Theorem \ref{yu-theorem-5-20-1}.]

      By Theorem
    \ref{yu-theorem-7-3-1}, we only need to show the sufficiency.
    For this purpose, we suppose that  there are
    constants $\beta>0$ and $C(\beta)>0$ that satisfy (\ref{yu-6-9-1}).
      We divide the remainder of the proof into two steps.

\vskip 5pt
   \noindent \emph{Step 1. We prove the stabilizability for system $[A,B]$ for case $(a)$.}

   First, $(a)$ in Theorem \ref{yu-theorem-5-20-1} implies that there is $C_A>0$
   such that
       $\|S^*(t)\|_{\mathcal{L}(H)}\leq C_A$ for each $t\geq 0$.

       Next, we arbitrarily fix $\tau>0$. Let
\begin{equation}\label{yu-5-18-20}
    \Theta_\tau(t):=
\begin{cases}
    \sin\left(\frac{\pi t}{\tau}\right),&\mbox{if}\;\;t\in[0,\tau],\\
    0,&\mbox{if}\;\;t\in\mathbb{R}\setminus[0,\tau],
\end{cases}
\end{equation}
    which satisfies
    \begin{equation}\label{3.13-7-23-w}
    \Theta_\tau(\cdot)\in H^1(\mathbb{R})
    \end{equation}
    and
\begin{equation}\label{yu-5-18-22}
    \Theta_{\tau}'(t)=
\begin{cases}
    \frac{\pi}{\tau}\cos\left(\frac{\pi t}{\tau}\right),&\mbox{if}\;\;t\in[0,\tau],\\
    0,&\mbox{if}\;\;t\in \mathbb{R}\setminus[0,\tau].
\end{cases}
\end{equation}
Now, we arbitrarily fix $\varphi\in H_1$ and define
\begin{equation}\label{yu-5-18-23}
  w(t):=\Theta_{\tau}(t)z(t),\;t\in \mathbb{R},\;\;\mbox{where}\;\;  z(t):=
    \begin{cases}
    S^*(t)\varphi,&\mbox{if}\;\;t\geq 0,\\
    0,&\mbox{if}\;\;t<0.
\end{cases}
\end{equation}
   By assumption, we have
       $\|z(t)\|_H\leq C_A\|\varphi\|_H$ for each $t\in\mathbb{R}$. The combination of this result,
       \eqref{yu-5-18-23} and \eqref{3.13-7-23-w} imply that
       $w\in H^1(\mathbb{R};H)$ and
    \begin{equation}\label{yu-5-18-24}
    w'(t)=
\begin{cases}
    A^*w(t)+\Theta'_{\tau}(t)z(t)&\mbox{if}\;\;t\geq 0,\\
    0,&\mbox{if}\;\;t<0.
\end{cases}
\end{equation}
             Thus, we can
             apply the Fourier transform to (\ref{yu-5-18-24}) to obtain
\begin{equation}\label{3.18-7-23w}
    (\mathrm{i}\varsigma I-A^*)\mathcal{F}[w](\varsigma)=\mathcal{F}[g](\varsigma)\;\;\mbox{a.e.}\;\;\varsigma\in\mathbb{R},
\end{equation}
    where $g(\cdot):=\Theta'_{\tau}(\cdot)z(\cdot)$, and  $\mathcal{F}[f]$ denotes the Fourier transform of $f\in L^2(\mathbb{R};H)$.
     Integrating (\ref{yu-6-9-1}) (where $\varphi$ and $\lambda$ are replaced by $\mathcal{F}[w](\varsigma)$ and $\mathrm{i}\varsigma$)  with respect to $\varsigma$ over $\mathbb{R}$, using
     \eqref{3.18-7-23w}, we obtain
\begin{equation}\label{yu-5-18-27}
    \int_{\mathbb{R}}\|\mathcal{F}[w](\varsigma)\|^2_Hd\varsigma\leq
    \beta^{-2}C(\beta)\left(
    \int_{\mathbb{R}}\|\mathcal{F}[g](\varsigma)\|^2_Hd\varsigma+
    \int_{\mathbb{R}}\|B^*\mathcal{F}[w](\varsigma)\|^2_Ud\varsigma\right).
\end{equation}
   It is clear that
   \begin{equation}\label{3.19WWW7-23}
   \mathcal{F}[g](\cdot)\in L^2(\mathbb{R};H).
   \end{equation}
    Thus, the first integral on the right-hand side of \eqref{yu-5-18-27} is  finite. We now claim
   \begin{equation}\label{3.20-7-23w}
   B^*\mathcal{F}[w](\cdot)\in L^2(\mathbb{R};U).
    \end{equation}
   Two facts ensuring \eqref{3.20-7-23w} are as follows: First,   we have
    \begin{equation}\label{3.21-7-23w}
   \mathcal{F}[B^*w](\cdot)\in L^2(\mathbb{R};U).
      \end{equation}
    Indeed, it follows from (\ref{3.13-7-23-w}) and assumption $(A_2)$ (see $(i)$ in Remark \ref{yu-b-remark-5-11-1})
    that
\begin{eqnarray*}\label{yu-5-18-28}
    \|B^*w(\cdot)\|_U&=&\|\Theta_{\tau}(\cdot)B^*(\rho_0I-A^*)^{-1}S^*(\cdot)(\rho_0I-A^*)\varphi\|_U
    \nonumber\\
    &\leq& C_{A} |\Theta_{\tau}(\cdot)|\|B^*(\rho_0I-A^*)^{-1}\|_{\mathcal{L}(H;U)}\|(\rho_0I-A^*)\varphi\|_H
    \in L^2(\mathbb{R}^+;\mathbb{R}),
\end{eqnarray*}
    while it follows from (\ref{yu-5-18-20}) and (\ref{yu-5-18-23})
     that $\|B^*w(t)\|_U=0$ when $t\in \overline{\mathbb{R}^-}$.
     These relations lead to \eqref{3.21-7-23w}. Second,
        one can check that for a.e. $\varsigma\in\mathbb{R}$,
\begin{eqnarray}\label{3.22-7-23w}
    \mathcal{F}[B^*w](\varsigma)&=&\int_{\mathbb{R}}e^{-\mathrm{i}\varsigma t}B^*w(t)dt=
    \int_{\mathbb{R}^+}e^{-\mathrm{i}\varsigma t}B^*(\rho_0I-A^*)^{-1}\Theta_{\tau}(t)S^*(t)(\rho_0I-A^*)\varphi dt\nonumber\\
    &=&B^*(\rho_0I-A^*)^{-1}\int_{\mathbb{R}^+}e^{-\mathrm{i}\varsigma t}\Theta_{\tau}(t)S^*(t)(\rho_0I-A^*)\varphi dt=B^*\int_{\mathbb{R}^+}e^{-\mathrm{i}\varsigma t}\Theta_{\tau}(t)S^*(t)\varphi dt\nonumber\\
    &=&B^*\int_{\mathbb{R}}e^{-\mathrm{i}\varsigma t}w(t)dt=B^*\mathcal{F}[w](\varsigma),
\end{eqnarray}
    where we used the fact $w(\cdot)\in L^2(\mathbb{R};H)$.
     Clearly, \eqref{3.20-7-23w} follows from \eqref{3.21-7-23w} and \eqref{3.22-7-23w}.

     Now, based on \eqref{3.19WWW7-23},  \eqref{3.20-7-23w}, (\ref{yu-5-18-27}), (\ref{yu-5-18-20}), (\ref{3.13-7-23-w}) and Plancherel's theorem, we obtain
\begin{equation}\label{yu-5-19-1}
    \int_0^{\tau}\|\Theta_{\tau}(t)z(t)\|_H^2\leq \beta^{-2}C(\beta)\left(\int_0^{\tau}\|B^*\Theta_{\tau}(t)z(t)\|_U^2dt
    +\int_0^{\tau}\|\Theta_{\tau}'(t)z(t)\|_H^2dt\right).
\end{equation}

We will use \eqref{yu-5-19-1} to obtain the weak observability, which leads to the stabilizability of $[A,B]$.
Indeed, since $\|S^*(t)\|_{\mathcal{L}(H)}\leq C_A$ for all $t\geq 0$, it follows from
    (\ref{yu-5-18-20}) that
\begin{eqnarray*}\label{yu-5-19-2}
    \frac{1}{2}\tau C_A^{-2}\|S^*(\tau)\varphi\|_H^2
    =C_A^{-2}\int_{\frac{\tau}{4}}^{\frac{3\tau}{4}}\|S^*(\tau)\varphi\|_H^2ds
    \leq \int_{\frac{\tau}{4}}^{\frac{3\tau}{4}}\|S^*(s)\varphi\|_H^2ds;
\end{eqnarray*}
   \begin{eqnarray*}\label{yu-5-19-3}
    \int_0^{\tau}\|\Theta_{\tau}(t)z(t)\|_H^2dt=\int_0^{\tau}\sin^2\left(\frac{\pi t}{\tau}
    \right)\|S^*(t)\varphi\|_H^2dt
    \geq \int_{\frac{\tau}{4}}^{\frac{3\tau}{4}}\sin^2\left(\frac{\pi t}{\tau}
    \right)\|S^*(t)\varphi\|_H^2dt\geq \frac{1}{2}\int_{\frac{\tau}{4}}^{\frac{3\tau}{4}}
    \|S^*(t)\varphi\|_H^2dt.
\end{eqnarray*}
    These relations yield that
\begin{eqnarray}\label{yu-5-19-4}
    \frac{1}{4}\tau C_A^{-2}\|S^*(\tau)\varphi\|_H^2\leq
    \int_0^\tau\|\Theta_{\tau}(t)z(t)\|^2_Hdt.
\end{eqnarray}
    Using (\ref{yu-5-18-20}) and (\ref{yu-5-18-22}), we further have
\begin{equation*}\label{yu-5-19-5}
    \int_0^{\tau}\|B^*\Theta_{\tau}(t)z(t)\|^2_Udt\leq
    \int_0^{\tau}\|B^*S^*(t)\varphi\|^2_Udt;
\end{equation*}
  \begin{equation*}\label{yu-5-19-6}
    \int_{0}^{\tau}\|\Theta_{\tau}'(t)z(t)\|_H^2dt
    \leq \tau^{-2}\pi^2\int_0^{\tau}\|S^*(t)\varphi\|_H^2dt\leq
    \tau^{-1}(\pi C_A)^{2}\|\varphi\|_H^2.
\end{equation*}
    With (\ref{yu-5-19-1}) and (\ref{yu-5-19-4}), these relations imply
\begin{equation}\label{yu-5-19-7}
    \|S^*(\tau)\varphi\|^2_H\leq 4\tau^{-1}\beta^{-2}
    C(\beta)C_A^2\int_0^\tau\|B^*S^*(t)\varphi\|_U^2dt
    +4\tau^{-2}\beta^{-2}C(\beta)\pi^2C_A^4\|\varphi\|_H^2.
\end{equation}
    Taking $\tilde{\tau}>0$ such that  $4\tilde{\tau}^{-2}\beta^{-2}C(\beta)\pi^2C_A^4\leq \frac{1}{2}$ in  (\ref{yu-5-19-7}) leads to
  \begin{equation*}\label{yu-5-19-8}
    \|S^*(\tilde{\tau})\varphi\|^2_H\leq 4\tilde{\tau}^{-1}\beta^{-2}
    C(\beta)C_A^2\int_0^{\tilde{\tau}}\|B^*S^*(t)\varphi\|_U^2dt
    +\frac{1}{2}\|\varphi\|_H^2.
\end{equation*}
    The combination of this result and conclusion $(i)$ in Lemma \ref{yu-bb-remark-6-21-1} imply  that
     system $[A,B]$ is  stabilizable.
\vskip 5pt
   \noindent \emph{Step 2. We prove the stabilizability for system $[A,B]$ for case $(b)$.}

   According to assumption $(b)$  (in Theorem \ref{yu-theorem-5-20-1}),  there are two closed subspaces $Q_1:=Q_1(\beta)$ and $Q_2:=Q_2(\beta)$ of $H$ satisfying  $(b_1)$-$(b_4)$ (where $\alpha$ is replaced by $\beta$). Based on $(b_1)$, we can define $P:H\to Q_1$ in the following manner:  $Pf=f_1$, for each $f\in H$,
   where $f=f_1+f_2$ with $f_j\in Q_j$ ($j=1,2$).
   Based on assumption $(b_2)$, one can directly check the following:
   \begin{equation}\label{3.27-7-24w}
   PS^*(t)=S^*(t)P,\; t\geq 0;\;\; PH_1\subset H_1;\;\; A^*P=PA^*\;\mbox{on}\;H_1.
   \end{equation}
   We write $S_1^*(\cdot):=PS^*(\cdot)$ and $S^*_2(\cdot):=(I-P)S^*(\cdot)$; $A_1^*=A^*|_{Q_1}$ and $A_2^*=A^*|_{Q_2}$. Based on \eqref{3.27-7-24w}, one can easily check that
   $S_j^*(\cdot) $ is the $C_0$-semigroup on $Q_j$, generated by $A^*_j$, $j=1,2$.

   Two facts are as follows: First, by $(b_3)$, we have that $A_1^*\in \mathcal{L}(Q_1)$,  which implies
   $D(A_1^*)=Q_1$. The combination of this result and assumption $(A_2)$ (see also $(i)$ in Remark \ref{yu-b-remark-5-11-1}) yield that
  \begin{equation}\label{3.28-7-24w}
  B^*P|_{Q_1}=(B^*(\rho_0I-A^*)^{-1})(\rho_0I-A^*_1)\in \mathcal{L}(Q_1;U).
 \end{equation}
   Second, by $(b_3)$, we have $\sigma(A_1^*)\subset \mathbb{C}_{-\beta}^+$ which implies
   \begin{equation}\label{3.29-7-24w}
   \overline{\mathbb{C}_{-\beta}^-}\subset \rho(A_1^*).
   \end{equation}
   Based on \eqref{3.29-7-24w} and  (\ref{yu-6-9-1}), for each $\lambda\in\mathbb{C}$, there is $C(\lambda)>0$ such that
\begin{equation}\label{3.30-7-24w}
    \|\varphi\|_H^2\leq C(\lambda)(\|(\lambda I-A_1^*)\varphi\|_H^2+\|B^*P\varphi\|_U^2)
    \;\;\mbox{for}\;\;\varphi\in Q_1.
\end{equation}
    (Indeed, if $\lambda\in\mathbb{C}_{-\beta}^+$, we can use (\ref{yu-6-9-1}) to find $C(\lambda)$ above, while when $\lambda\notin\mathbb{C}_{-\beta}^+$, we can choose $C(\lambda)=\|(\lambda I-A_1^*)^{-1}\|_{\mathcal{L}(H)}^2$ because of \eqref{3.29-7-24w}.)

    Now, based on \eqref{3.30-7-24w}, \eqref{3.28-7-24w} and boundedness of $A_1^*$, we can apply
     Lemma \ref{yu-lemma-7-6-1} to find $T>0$ and $C(T)>0$ such that
\begin{equation*}\label{yu-7-12-10}
    \|S^*_1(T)\varphi\|_H^2\leq C(T)\int_0^T\|B^*PS^*_1(t)\varphi\|_U^2dt
    \;\;\mbox{for}\;\;\varphi\in Q_1.
\end{equation*}
    The combination of this result and $(b_2)$ yield
    \begin{eqnarray}\label{yu-5-26-9}
    &&|PS^*(T)\varphi\|^2_H=\|S_1^*(T)P\varphi\|_H^2\nonumber\\
    &&\leq C(T)\int_0^T\|B^*PS_1^*(t)P\varphi\|_U^2dt
    =C(T)\int_0^T\|B^*PS^*(t)\varphi\|_U^2dt\;\;\mbox{for}\;\;\varphi\in H_1.
\end{eqnarray}
    Based on $(b_4)$ and (\ref{yu-5-26-9}),  using a similar method as in \emph{Sub-step 1.2} of the proof of  Theorem
    \ref{yu-theorem-6-30-1}, there are $\widehat{T}>0$ and $C(\widehat{T})>0$ such that
\begin{equation*}\label{yu-5-29-2}
    \|S^*(\widehat{T})\varphi\|_H^2\leq \widehat{C}(\widehat{T})\int_0^{\widehat{T}}\|B^*S^*(t)\|_U^2dt+\frac{1}{2}
    \|\varphi\|_H^2\;\;\mbox{for}\;\;\varphi\in H_1,
\end{equation*}
     The combination of this result and $(i)$ of Lemma \ref{yu-bb-remark-6-21-1} lead to the stabilizability of system $[A,B]$.
\par
    Hence, we have completed the proof of Theorem \ref{yu-theorem-5-20-1}.
\end{proof}

\section{Applications}
This section provides several applications of our main theorems to specific control PDEs. We start with introducing the concept of `thick sets': {\it We say a measurable subset $E\subset \mathbb{R}^N$ (with $N\in\mathbb{N}^+$) to be thick,
if  there is $\varepsilon>0$ and $L>0$ such that
\begin{equation*}\label{yu-6-5-12}
    |E\cap Q_L(x)|\geq \varepsilon L^N\;\;\mbox{for each}\;\;x\in\mathbb{R}^N,
\end{equation*}
    where $Q_L(x)$ denotes the closed cube in $\mathbb{R}^N$, centered at $x$ and of side length $L$,
    and $|E\cap Q_L(x)|$ denotes the Lebesgue measure of $E\cap Q_L(x)$. }
    Then, we quote the following lemma, which is related to the thick sets and will be used later:
 \begin{lemma}\label{yu-lemma-6-15-1}
    (\cite[Theorem 1]{Kovrijkine}) If $\omega$ is a thick set, then for each $R>0$, there exists  $C(R,\omega)>0$ such that, for each $f\in L^2(\mathbb{R}^N)$ with $\mbox{supp}(\mathcal{F}[f])\subset [-R,R]^N$, the following estimate holds
\begin{equation*}\label{yu-6-15-3}
    \|f\|_{L^2(\mathbb{R}^N)}\leq C(R,\omega)\|\chi_{\omega}f\|_{L^2(\mathbb{R}^N)}.
\end{equation*}
\end{lemma}
%
%


\subsection{Ginzburg--Landau equation in $\mathbb{R}^N$}\label{yu-sec-4-2}
    Let $a\in\mathbb{R}^+$ and $b\in\mathbb{R}$.
     Let $\omega\subset \mathbb{R}^N$ (with $N\in\mathbb{N}^+$) be a  measurable set with its characteristic function $\chi_{\omega}$.
          We consider the controlled Ginzburg--Landau equation in $\mathbb{R}^N$:
\begin{equation}\label{yu-6-15-1}
\begin{cases}
    y_t=(a+\mathrm{i}b)\triangle y+\chi_{\omega}u&\mbox{in}\;\;\mathbb{R}^+\times \mathbb{R}^N,\\
    y(0,\cdot)=y_0(\cdot)\in L^2(\mathbb{R}^N),
\end{cases}
\end{equation}
    where  $u\in L^2(\mathbb{R}^+;L^2(\mathbb{R}^N))$.
Equation \eqref{yu-6-15-1} can be put into our framework by setting:
   $H=U=L^2(\mathbb{R}^N)$;  $A:=(a+\mathrm{i}b)\triangle$, with its domain $H^2(\mathbb{R}^N)$;
     $B:=\chi_{\omega}$. One can directly check that  $(A_1)$-$(A_3)$ are true.
      One can also check that $A$ is normal (using Fourier transform) and
       generates an analytic semigroup.
       The latter and $(ii)$ of Remark \ref{remark 1.10-7-24w} yield that $\sigma(A)\cap \mathbb{C}_{-\gamma}^+$ is bounded for each $\gamma>0$. Therefore, Theorem \ref{yu-theorem-6-30-1}
       can be applied. It provides the following results:

\begin{theorem}\label{yu-theorem-6-15-1}
    If $\omega$ is a thick set, then equation \eqref{yu-6-15-1} is rapidly stabilizable.
\end{theorem}
\begin{proof}
   According to $(ii)$ in Theorem \ref{yu-theorem-6-30-1} and  Definition \ref{yu-definition-6-9-1}, we only need to show the following: For each $\beta>0$, there is $C(\beta)>0$ such that
\begin{equation}\label{yu-b-6-15-1}
    \|\varphi\|^2_{L^2(\mathbb{R}^N)}\leq C(\beta)(\|(\lambda I-(a-\mathrm{i}b)\triangle)\varphi\|^2_{L^2(\mathbb{R}^N)}
    +\|\chi_{\omega}\varphi\|_{L^2(\mathbb{R}^N)}^2)
    \;\;\mbox{for}\;\;\lambda\in\mathbb{C}_{-\beta}^+,\;\varphi\in H^2(\mathbb{R}^N).
\end{equation}
    For this purpose, we arbitrarily fix  $\varphi\in H^2(\mathbb{R}^N)$, $\beta>0$ and $\lambda\in\mathbb{C}_{-\beta}^+$. Plancherel's theorem implies that
\begin{eqnarray*}\label{yu-6-15-10}
    \|(\lambda I-(a-\mathrm{i}b)\triangle)\varphi\|^2_{L^2(\mathbb{R}^N)}
    &=&\|(\lambda+(a-\mathrm{i}b)|\xi|^2)\mathcal{F}[\varphi]\|_{L^2(\mathbb{R}^N)}^2\nonumber\\
    &=&\int_{\mathbb{R}^N}
    \left(|\mbox{Re}\lambda+a|\xi|^2|^2+|\mbox{Im}\lambda-b|\xi|^2|^2\right)
    |\mathcal{F}[\varphi](\xi)|^2d\xi\nonumber\\
    &\geq&\beta^2\int_{\mathbb{R}^N}\chi_{|\xi|\geq \sqrt{2a^{-1}\beta}}
    |\mathcal{F}[\varphi](\xi)|^2d\xi.
\end{eqnarray*}
     The combination of this result, Lemma \ref{yu-lemma-6-15-1}, and the Plancherel theorem yield that there is $C_1(\beta,a,\omega)>0$ such that
\begin{eqnarray*}\label{yu-7-4-1}
    \|\varphi\|^2_{L^2(\mathbb{R}^N)}
    &\leq& \beta^{-2}\|(\lambda I-(a-\mathrm{i}b)\triangle)\varphi\|^2_{L^2(\mathbb{R}^N)}
    +\int_{\mathbb{R}^N}\chi_{|\xi|\leq\sqrt{2a^{-1}\beta}}
    |\mathcal{F}[\varphi](\xi)|^2d\xi\nonumber\\
    &\leq&\beta^{-2}\|(\lambda I-(a-\mathrm{i}b)\triangle)\varphi\|^2_{L^2(\mathbb{R}^N)}
    +C_1(\beta,a,\omega)\int_{\mathbb{R}^N}\chi_{\omega}(x)|
    \mathcal{F}^{-1}[\chi_{|\xi|\leq\sqrt{2a^{-1}\beta}}
    \mathcal{F}[\varphi]](x)|^2dx\nonumber\\
    &\leq&\beta^{-2}\|(\lambda I-(a-\mathrm{i}b)\triangle)\varphi\|^2_{L^2(\mathbb{R}^N)}
    +2C_1(\beta,a,\omega)\int_{\mathbb{R}^N}\chi_{\omega}(x)|\varphi(x)|^2dx\nonumber\\
    &\;&+2C_1(\beta,a,\omega)\int_{\mathbb{R}^N}|
    \mathcal{F}^{-1}[\chi_{|\xi|\geq\sqrt{2a^{-1}\beta}}
    |\mathcal{F}[\varphi]](x)|^2dx\nonumber\\
    &\leq&\beta^{-2}(2C_1(\beta,a,\omega)+1)\|(\lambda I-(a-\mathrm{i}b)\triangle)\varphi\|^2_{L^2(\mathbb{R}^N)}
    +2C_1(\beta,a,\omega)\|\chi_{\omega}\varphi\|^2_{L^2(\mathbb{R}^N)},
\end{eqnarray*}
    which leads to
    (\ref{yu-b-6-15-1}) and completes the proof.
\end{proof}

\subsection{Fractional heat equation in $\mathbb{R}^N$}\label{yu-7-7-1-sc}
Let $s\in(0,1)$ and $N\in\mathbb{N}^+$. Let $\omega\subset \mathbb{R}^N$ be a  measurable set with its characteristic function $\chi_{\omega}$.
 We consider the following controlled fractional heat equation:
\begin{equation}\label{yu-6-5-10}
\begin{cases}
    \partial_ty+(-\triangle)^{\frac{s}{2}}y=\chi_\omega u&\mbox{in}\;\;\mathbb{R}^+\times \mathbb{R}^N,\\
    y(0,\cdot)=y_0(\cdot)\in L^2(\mathbb{R}^N),
\end{cases}
\end{equation}
    where
    $u\in L^2(\mathbb{R}^+;\mathbb{R}^N)$ and $(-\triangle)^{\frac{s}{2}}$ is defined by
\begin{equation}\label{yu-6-5-11}
    (-\triangle)^{\frac{s}{2}}f:=
    \mathcal{F}^{-1}[|\xi|^s\mathcal{F}[f]], \;\;f\in C_c^{\infty}(\mathbb{R}^N).
\end{equation}
Equation \eqref{yu-6-5-10} can be put into our framework by setting:
  $H=U=L^2(\mathbb{R}^N)$;   $A:=-(-\triangle)^{\frac{s}{2}}$ with domain $D(A):=H^s(\mathbb{R}^N)$;
   $B:=\chi_\omega(\in\mathcal{L}(L^2(\mathbb{R}^N)))$. One can easily check  that  $(A_1)$-$(A_3)$ hold, $A^*=A$, and $A$ generates an analytic semigroup. The latter and $(ii)$
   of Remark \ref{remark 1.10-7-24w} yield that $\sigma(A)\cap \mathbb{C}_{-\gamma}^+$ is bounded for each $\gamma>0$. Therefore, Theorem \ref{yu-theorem-6-30-1}
       can be applied. Moreover, the spectral measure $ E^{A^*}$ (corresponding to $A^*$) is provided as follows: for
            each Borel set $\Omega\subset\mathbb{C}$,
                        \begin{equation}\label{yu-7-3-10}
   [E^{A^*}(\Omega)]f:=\mathcal{F}^{-1}[\chi_{-|\xi|^s\in (\Omega\cap \overline{\mathbb{R}^-})}\mathcal{F}[f]]
    \;\;\mbox{for}\;\;f\in L^2(\mathbb{R}^N).
\end{equation}

Regarding equation \eqref{yu-6-5-10}, we have the following: {\it Fact 1:}
If $\omega$ is thick, then
equation \eqref{yu-6-5-10} is not null controllable in general (see \cite{Koenig}).
{\it Fact 2:}
Equation \eqref{yu-6-5-10} is rapidly stabilizable if and only if $\omega$ is thick (see \cite[Theorem 4.5]{Liu-Wang-Xu-Yu}), where it was proven by the weak observability inequality.
Next, we utilize $(ii)$ of Theorem \ref{yu-theorem-6-30-1} to provide a proof for the sufficiency by the frequency-domain inequality.
%
\begin{proposition}
    If $\omega$ is a thick set, then for each $\beta>0$, there is $C(\beta)>0$ such that
\begin{equation}\label{yu-6-6-00}
    \|\varphi\|_{L^2(\mathbb{R}^N)}^2\leq \frac{C(\beta)}{(\mbox{Re}\;\lambda+\beta)^2}
    \left(\|(\lambda I -A^*)\varphi\|^2_{L^2(\mathbb{R}^N)}+
    \|B^*\varphi\|^2_{L^2(\mathbb{R}^N)}\right)
    \;\;\mbox{for}\;\;\lambda\in\mathbb{C}_{-\beta}^+,\;\varphi\in H^s(\mathbb{R}^N).
\end{equation}
\end{proposition}
\begin{proof}
    We arbitrarily fix $\beta>0$. Let $k:=k(\beta)=[\beta]+2$. We set
    $E_k:=E^{A^*}({\{z\in\mathbb{C}:\mbox{Re}z\in[-k,0]\}})$, where $E^{A^*}$ is provided in (\ref{yu-7-3-10}).
 Then, \cite[Lemma 3.1]{Huang-Wang} implies that
\begin{equation}\label{yu-6-5-21}
    \|E_k\varphi\|^2_{L^2(\mathbb{R}^N)}\leq e^{C_0k^{\frac{1}{s}}}\|B^*E_k\varphi\|^2_{L^2(\mathbb{R}^N)}
    \;\;\mbox{for}\;\;\varphi\in L^2(\mathbb{R}^N),
\end{equation}
    where $C_0>0$ is a constant, which is independent of $k$. Let $A^*_k:=A^*+(k-\frac{1}{2})I$.
      Now we claim
\begin{equation}\label{yu-6-5-22}
    \|(I-E_k)\varphi\|_{L^2(\mathbb{R}^N)}^2
    \leq \frac{1}{(\mbox{Re}\lambda+\frac{1}{2})^2}\|(\lambda I-A^*_k)(I-E_k)\varphi\|_{L^2(\mathbb{R}^N)}
    \;\;\mbox{for}\;\;\lambda\in\mathbb{C}_{-\frac{1}{2}}^+,\;\varphi\in H^s(\mathbb{R}^N).
\end{equation}
   For this purpose, we arbitrarily fix $\lambda\in\mathbb{C}_{-\frac{1}{2}}^+$, $\varphi\in H^s(\mathbb{R}^N)$. Since  $A=A^*$, it follows from (\ref{yu-6-5-11}) and (\ref{yu-7-3-10}) that
\begin{equation*}\label{yu-6-5-23}
    (\lambda I-A^*_k)(I-E_k)\varphi= \mathcal{F}^{-1}\Big[\chi_{\{|\xi|^s>k\}}\Big(\lambda
    -k+\frac{1}{2}+|\xi|^s\Big)\mathcal{F}[\varphi]\Big].
\end{equation*}
    The combination of this result and the Plancherel theorem yield
\begin{eqnarray*}\label{yu-6-5-24}
    &\;&\|(\lambda I-A^*_k)(I-E_k)\varphi\|_{L^2(\mathbb{R}^N)}^2=
    \int_{\mathbb{R}^N}\Big|\chi_{\{|\xi|^s>k\}}\Big(\lambda-k+\frac{1}{2}+|\xi|^s\Big)
    \mathcal{F}[\varphi](\xi)
    \Big|^2d\xi\nonumber\\
    &\geq&\Big|\lambda+\frac{1}{2}\Big|^2
    \int_{\mathbb{R}^N}|\chi_{\{|\xi|^s>k\}}\mathcal{F}[\varphi]
    (\xi)|^2d\xi\geq\Big(\mbox{Re}\lambda+\frac{1}{2}\Big)^2\|(I-E_k)\varphi\|^2_{L^2(\mathbb{R}^N)},
    \end{eqnarray*}
    which leads to  (\ref{yu-6-5-22}).

    Next, since $A^*E_k=E_kA^*$, (\ref{yu-6-5-21}) and (\ref{yu-6-5-22}) imply that for $\lambda\in\mathbb{C}_{-\frac{1}{4}}^+$ and $\varphi\in H^s(\mathbb{R}^N)$,
\begin{eqnarray*}\label{yu-6-5-25}
    \|\varphi\|_{L^2(\mathbb{R}^N)}^2&\leq& \|(I-E_k)\varphi\|_{L^2(\mathbb{R}^N)}^2+\|E_k\varphi\|_{L^2(\mathbb{R}^N)}^2
    \leq \|(I-E_k)\varphi\|_{L^2(\mathbb{R}^N)}^2+
    e^{C_0k^{\frac{1}{s}}}\|B^*E_k\varphi\|^2_{L^2(\mathbb{R}^N)}
    \nonumber\\
    &\leq& \Big(1+2e^{C_0k^{\frac{1}{s}}}
    \Big)\|(I-E_k)\varphi\|^2_{L^2(\mathbb{R}^N)}
    +2e^{C_0k^{\frac{1}{s}}}\|B^*\varphi\|^2_{L^2(\mathbb{R}^N)}\nonumber\\
    &\leq& 16\Big(1+2e^{C_0k^{\frac{1}{s}}}
    \Big)\|(\lambda I-A^*_k)\varphi\|^2_{L^2(\mathbb{R}^N)}
    +2e^{C_0k^{\frac{1}{s}}}\|B^*\varphi\|^2_{L^2(\mathbb{R}^N)}.
\end{eqnarray*}
          The combination of this result and Definition \ref{yu-definition-6-9-1} (see also   Proposition \ref{yu-prop-6-3-1} in Appendix)  imply that
     there are $\gamma>0$ and $C(\gamma,k)>0$ such that
\begin{equation}\label{yu-6-6-1}
    \|\varphi\|_{L^2(\mathbb{R}^N)}^2\leq \frac{C(\gamma,k)}{(\mbox{Re}\,\eta+\gamma)^2}
    \left(\|(\eta I-A^*_k)\varphi\|^2_{L^2(\mathbb{R}^N)}+
    \|B^*\varphi\|^2_{L^2(\mathbb{R}^N)}\right)
    \;\;\mbox{for}\;\;\eta\in\mathbb{C}_{-\gamma}^+,\;\varphi\in H^s(\mathbb{R}^N).
\end{equation}
    Since $\beta+1\leq k$,  letting $\eta=\lambda+k-\frac{1}{2}$ in (\ref{yu-6-6-1}) leads to (\ref{yu-6-6-00}) with $C(\beta):=C(\gamma,k)$. This completes the proof.
\end{proof}

\subsection{One-dimensional heat equation with point-wise controls}\label{yu-sec-4-4}
  Let  $c>\pi^2$, $x_0\in(0,1)$, and $\delta(\cdot)$ be the Dirac function at $x=0\in \mathbb{R}$. We consider the following  heat equation with point-wise controls:
\begin{equation}\label{yu-6-19-7}
\begin{cases}
    y_t=(\partial_x^2+c) y+\delta(\cdot-x_0)u &\mbox{in}\;\;\mathbb{R}^+\times (0,1), \\
    y(\cdot,0)=y(\cdot,1)=0&\mbox{in}\;\;\mathbb{R}^+,\\
    y(0,\cdot)=y_0(\cdot)\in L^2(0,1),
\end{cases}
\end{equation}
    where  $u\in L^2(\mathbb{R}^+)$. Equation \eqref{yu-6-19-7} can be put into our framework by setting:
     $H:=L^2(0,1)$; $U:=\mathbb{R}$; $A:=\partial_x^2+c$ (with its domain $D(A):=H_1=H_0^1(0,1)\cap H^2(0,1)$);
    $B:=\delta(\cdot-x_0)$  (which is considered as a bounded operator from $\mathbb{R}$ to in $H_{-1}$, see $(i)$ in Remark \ref{yu-b-remark-5-11-1} for the definition of $H_{-1}$).
    One can easily check that assumptions $(A_1)$-$(A_3)$ are true, $A$ is self-adjoint, and $A$ generates an analytic semigroup
    (see \cite[Example 4.3]{Liu-Wang-Xu-Yu}). The latter and $(ii)$
   of Remark \ref{remark 1.10-7-24w} yield that $\sigma(A)\cap \mathbb{C}_{-\gamma}^+$ is bounded for each $\gamma>0$. Therefore, Theorem \ref{yu-theorem-6-30-1}
       can be applied. Moreover, since $A$ generates a compact semigroup, assumption $(b)$ in Theorem \ref{yu-theorem-5-20-1} holds (see $(ii)$ in Remark \ref{yu-remark-6-9-1}). Thus,  Theorem \ref{yu-theorem-5-20-1} can also be applied.

    Regarding equation \eqref{yu-6-19-7},
    we have the following: {\it Fact 1:} For some
    irrational number $x_0\in (0,1)$,
 equation \eqref{yu-6-19-7} is not null controllable (see \cite[Example 4.3]{Liu-Wang-Xu-Yu}).
    {\it Fact 2:}  Equation \eqref{yu-6-19-7} is rapidly stabilizable if and only if $x_0\in (0,1)$ is irrational (see \cite[Theorem 4.9]{Liu-Wang-Xu-Yu}). Thus, a natural question is {\it what happens about  equation \eqref{yu-6-19-7} when
    $x_0$ is rational?}
     We will give  the answer for this question by using Theorem \ref{yu-theorem-6-30-1} or Theorem \ref{yu-theorem-5-20-1}.

\begin{theorem}\label{yu-theorem-7-13-1}
    System \eqref{yu-6-19-7} is stabilizable in $L^2(0,1)$ if and only if $x_0\notin\{k/n\in(0,1):k\in\mathbb{N}^+, n=1,2,\ldots,[\sqrt{c}/\pi]\}$.
\end{theorem}
\begin{proof}
 First of all, we give two facts. Fact One: The eigenvalues and the corresponding normalized eigenfunctions of $A^*$ are as follows:
$\lambda_n:=-(n\pi)^2+c$; $e_n(x):=\sqrt{2}\sin(n\pi x)$, $x\in (0,1)$, $n\in \mathbb{N}^+$. Fact Two:
 $B^*\varphi=\varphi(x_0)$ for $\varphi\in H_1$ which is bounded from $H_1$ to $\mathbb{R}$.

 Now we prove the necessity.
By contradiction, we suppose that there is $x_0=k/n\in (0,1)$ (with $n\in\{1,2,\ldots,[\sqrt{c}/\pi]\}$, $k\in \mathbb{N}^+$) such that  system \eqref{yu-6-19-7} is stabilizable. Then, we have $e_n(x_0)=\sqrt{2}\sin(n\pi x_0)=0$,
which implies that $B^*e_n=0$. Since $1\leq n\leq [\sqrt{c}/\pi]$, one can directly check that
\begin{eqnarray*}
    \lambda_n=-(n\pi)^2+c\geq -([\sqrt{c}/\pi]\pi)^2+c\geq 0,
    \;\;\mbox{and thus}\;\;\lambda_n\in\mathbb{C}_0^+.
   \end{eqnarray*}
    For any $\beta>0$, the right-hand side of (\ref{yu-6-9-1}) (where $\lambda=\lambda_n$ and
    $\varphi=e_n$) is $0$, while the left-hand side of (\ref{yu-6-9-1}) (where $\lambda=\lambda_n$ and
    $\varphi=e_n$) is $1$. So (\ref{yu-6-9-1}) is not true for any $\beta>0$. Thus, it follows from
     $(i)$ of Theorem \ref{yu-theorem-6-30-1} (or Theorem \ref{yu-theorem-5-20-1}) that system \eqref{yu-6-19-7} is not  stabilizable, which causes a contradiction and completes the proof of  the necessity.

   We next prove the sufficiency.  We will show that if $x_0\notin\{k/n\in(0,1):k\in\mathbb{N}^+, n=1,2,\ldots,[\sqrt{c}/\pi]\}$, then there is $\beta>0$ such that the pair $(A^*,B^*)$ satisfies
   (HESI)$_\beta$. We let  $n^*:=[\sqrt{c}/\pi]$. Since $x_0\notin\{k/n\in(0,1):k\in\mathbb{N}^+, n=1,2,\ldots,n^*\}$, we have
\begin{equation}\label{yu-7-13-20}
    e_n(x_0)\neq 0\;\;\mbox{for each}\;\;n\in\{1,2,\ldots, n^*\}.
\end{equation}
    We define the following projection operator:
    $P_{n^*}\varphi:=\sum_{n=1}^{n^*}\langle \varphi, e_n\rangle_{H} e_n$ ($\varphi\in H)$.
    Let
   $\beta:=-\frac{1}{2}\lambda_{n^*+1}$.
    Since  $n^*=[\sqrt{c}/\pi]$, one can directly check that  $\beta>0$.
     Now, we claim that there is $C(\beta)>0$ such that
\begin{equation}\label{yu-7-13-28}
    \|\varphi\|_{H}^2\leq C(\beta)(\|(\lambda I-A^*)\varphi\|^2_{H}+\|B^*\varphi\|_U^2)
    \;\;\mbox{for}\;\;\lambda\in \mathbb{C}_{-\beta}^+,\;\varphi\in H_1.
\end{equation}
    Indeed, based on (\ref{yu-7-13-20}), one can  easily check that
    $\mbox{Ker}\left((\lambda I-A^*P_{n^*},B^*P_{n^*})^\top|_{P_{n^*}H}\right)=\{0\}
    \;\;\mbox{for}\;\;\lambda\in\mathbb{C}$.
    Since $P_{n^*}H$ and $U$ are finite-dimensional, this, along with Kalman's controllability condition, yields
    that system
    $[P_{n^*}A,P_{n^*}B]$ is controllable (consequently, is rapidly stabilizable)  on $P_{n^*}H$.
    Thus, by $(ii)$ in Theorem \ref{yu-theorem-7-3-1}, there is $C_0(\beta)>0$ such that
    \begin{equation}\label{yu-7-13-23}
    \|P_{n^*}\varphi\|^2_{H}\leq C_0(\beta)(\|(\lambda I-A^*)P_{n^*}\varphi\|_{H}^2
    +\|B^*P_{n^*}\varphi\|_U^2)\;\;\mbox{for}\;\;\lambda\in\mathbb{C}_{-\beta}^+,\;\;\varphi\in H_1.
\end{equation}
    Meanwhile, by the definition of $\beta$,  one can directly check that, for $\lambda\in\mathbb{C}_{-\beta}^+$,
\begin{equation}\label{yu-12-16-1}
    \|(\lambda I-A^*)\varphi\|_{H}^2\geq \|(\lambda I-A^*)(I-P_{n^*})\varphi\|_{H}^2
    \geq \beta^2\|(I-P_{n^*})\varphi\|_H\;\;\mbox{for}\;\;\varphi\in H_1.
\end{equation}
    Moreover, by the H\"{o}lder inequality, there exists   $C_1(\beta)>0$ such that, for $\lambda\in \mathbb{C}_{-\beta}^+$
    and $\varphi\in H_1$,
\begin{eqnarray*}\label{yu-7-14-6}
    \|B^*P_{n^*}\varphi\|_U^2&\leq& 2\|B^*\varphi\|_U^2+2\|B^*(I-P_{n^*})\varphi\|_U^2
    \leq 2\|B^*\varphi\|_U^2+4\Big|\sum_{n=n^*+1}^{+\infty}a_n\sin (n\pi x_0)\Big|^2\nonumber\\
    &\leq&2\|B^*\varphi\|_U^2+4\Big(\sum_{n=n^*+1}^{+\infty}|a_n(\lambda-\lambda_n)|^2\Big)
    \Big(\sum_{n=n^*+1}^{+\infty}|\beta+\lambda_n|^{-2}\Big)\nonumber\\
    &\leq& 2\|B^*\varphi\|_U^2+C_1(\beta)\|(\lambda I-A^*)\varphi\|_{H}^2.
\end{eqnarray*}
   With (\ref{yu-7-13-23}) and (\ref{yu-12-16-1}), the above yields that for any $\lambda\in \mathbb{C}_{-\beta}^+$
    and $\varphi\in H_1$,
\begin{eqnarray*}
    \|\varphi\|_{H}^2\leq \|P_{n^*}\varphi\|_H^2+\|(I-P_{n^*})\varphi\|_H^2
    \leq (1+C_0(\beta))(1+\beta^{-2}+C_1(\beta))\|(\lambda I-A^*)\varphi\|_{H}^2
    +2C_0(\beta)\|B^*\varphi\|_U^2,
\end{eqnarray*}
which leads to (\ref{yu-7-13-28}). So  $(A^*,B^*)$ satisfies (HESI)$_\beta$.
Then, according to
     $(i)$ in Theorem \ref{yu-theorem-6-30-1} (or Theorem \ref{yu-theorem-5-20-1}), system (\ref{yu-6-19-7}) is  stabilizable. Thus,  the sufficiency has been proven.

\par
    Hence, we have completed the proof of Theorem \ref{yu-theorem-7-13-1}.
\end{proof}
\subsection{Heat equation with time delay in $\mathbb{R}^N$}\label{sect-55}
 Let   $\tau>0$, $N\in\mathbb{N}$. Let  $\omega\subset \mathbb{R}^N$ be a  subset of positive measure.
  Let $\chi_{\omega}$ be the  characteristic function of $\omega$.
    We consider the following
   controlled heat equation with time delay in $\mathbb{R}^N$:
\begin{equation}\label{yu-12-7-1}
\begin{cases}
    y_t(t,x)=(\triangle-I) y(t,x)+y(t-\tau,x)+\chi_{\omega}u(t,x),&(t,x)\in\mathbb{R}^+\times\mathbb{R}^N,\\
    y(0,x)=y_0(x),&x\in\mathbb{R}^N,\\
    y(t,x)=f(t,x),&(t,x)\in (-\tau,0)\times\mathbb{R}^N,
\end{cases}
\end{equation}
    where $I$ is the identity operator on $L^2(\mathbb{R}^N)$,
     $y_0\in L^2(\mathbb{R}^N)$,
    $f\in L^2(-\tau,0;L^2(\mathbb{R}^N))$ and  $u\in L^2(\mathbb{R}^+;L^2(\mathbb{R}^N))$.
   It is well known that (\ref{yu-12-7-1}) has a unique solution in $C([0,+\infty);L^2(\mathbb{R}^N))$ (see \cite[Theorem 2.1]{Nakagiri-1988}).
  We denote this solution by $y(\cdot,\cdot;y_0,f,u)$ if it is viewed as a real-valued function of $(t,x)$, while by $y(\cdot;y_0,f,u)$
  if it is viewed as an $L^2(\mathbb{R}^N)$-valued function of $t$.

\par
    In this subsection, we will apply Theorem \ref{yu-theorem-5-20-1} to show the stabilizability
    of (\ref{yu-12-7-1}) in the  sense of the following definition: \color{black}
\begin{definition}\label{yu-def-12-18-1}
    System (\ref{yu-12-7-1}) is said to be stabilizable if there is $\alpha>0$, $C>0$ and a feedback law $K\in\mathcal{L}(L^2(\mathbb{R}^N)
    \times L^2(-\tau,0;L^2(\mathbb{R}^N));L^2(\mathbb{R}^N))$ such that  for any $y_0\in L^2(\mathbb{R}^N)$ and $f\in L^2(-\tau,0;L^2(\mathbb{R}^N))$, the solution $y_K(\cdot;y_0,f)$ of the closed-loop system:
\begin{equation}\label{yu-12-20-100}
\begin{cases}
    y_t(t,x)=(\triangle-I)y(t,x)+y(t-\tau,x)+\chi_{\omega}[K(y(t),
    y(t+\cdot))^\top](x),&(t,x)\in\mathbb{R}^+
    \times \mathbb{R}^N,\\
    y(0,x)=y_0(x),&x\in\mathbb{R}^N,\\
    y(t,x)=f(t,x),&(t,x)\in (-\tau,0)\times\mathbb{R}^N,
\end{cases}
\end{equation}
    satisfies that
\begin{equation}\label{yu-12-20-101}
    \|y_K(t;y_0,f)\|_{L^2(\mathbb{R}^N)}\leq C e^{-\alpha t}(\|y_0\|_{L^2(\mathbb{R}^N)}
    +\|f\|_{L^2(-\tau,0;L^2(\mathbb{R}^N))})\;\;\mbox{for}\;\;t\in\mathbb{R}^+.
\end{equation}
\end{definition}
   For this purpose, we first embed system (\ref{yu-12-7-1})  into our  framework.  Let $\mathcal{U}:=L^2(\mathbb{R}^N)$ and $\mathcal{H}:=L^2(\mathbb{R}^N)\times L^2(0,1;L^2(\mathbb{R}^N))$ with the inner product:
$$
    \langle(\varphi_1,\varphi_2)^\top,(f_1,f_2)^\top\rangle_{\mathcal{H}}
    :=\langle\varphi_1,f_1\rangle_{L^2(\mathbb{R}^N)}+\tau\int_0^1\langle \varphi_2(\cdot,\rho)
    f_2(\cdot,\rho)\rangle_{L^2(\mathbb{R}^N)}d\rho
$$
    for $(\varphi_1,\varphi_2)^\top$, $(f_1,f_2)^\top
    \in L^2(\mathbb{R}^N)\times L^2(0,1;L^2(\mathbb{R}^N))$.
   We define  an unbounded operator on $\mathcal{H}$ by
\begin{equation}\label{yu-12-7-4-b}
\begin{cases}
    \mathcal{A}\varphi:=\left(
                          \begin{array}{c}
                            \mathcal{A}_1\varphi \\
                            \mathcal{A}_2\varphi \\
                          \end{array}
                        \right)\;\;\mbox{for}\;\;\varphi=\left(
                                                  \begin{array}{c}
                                                    \varphi_1 \\
                                                    \varphi_2 \\
                                                  \end{array}
                                                \right)
                      \in D(\mathcal{A}),\\
                        \;\;\mbox{with}\;\;\left(
                          \begin{array}{c}
                            [\mathcal{A}_1\varphi](x) \\
                            [\mathcal{A}_2\varphi](x,\rho) \\
                          \end{array}
                        \right)=
               \left(
                         \begin{array}{c}
                           (\triangle-I)\varphi_1(x)+\varphi_2(x,1) \\
                            -\tau^{-1}\partial_\rho\varphi_2(x,\rho)\\
                         \end{array}
                       \right),\;\;(x,\rho)\in\mathbb{R}^N\times(0,1),\\
                       D(\mathcal{A}):=
    \{(f_1,f_2)^{\top}\in \mathcal{H}: f_1\in  H^2(\mathbb{R}^N),
    \;f_2\in H^1(0,1;L^2(\mathbb{R}^N)),
    \;f_1(x)=f_2(x,0)\;\mbox{a.e.}\;x\in\mathbb{R}^N\}.
\end{cases}
\end{equation}
 Let  $\mathcal{B}:=(\chi_{\omega},0)^{\top}$.
  Then one can easily check that  $D(\mathcal{A})$ is dense in $\mathcal{H}$ and $\mathcal{B}\in \mathcal{L}(\mathcal{U};\mathcal{H})$.
 Moreover, we have the following result:
 \begin{lemma}\label{yu-lemma-12-8-1}
   The operator $\mathcal{A}$  generates a $C_0$-semigroup of contractions on
   $\mathcal{H}$.
\end{lemma}
 \color{black}
   \begin{proof} We first prove that  $\mathcal{A}$ is dissipative, i.e.,
\begin{equation}\label{yu-12-7-10-b-b}
    \mbox{Re}\langle \mathcal{A}\varphi,\varphi\rangle_{\mathcal{H}}\leq 0
    \;\;\mbox{for}\;\;\varphi\in D(\mathcal{A}).
\end{equation}
    Indeed,  it follows from  (\ref{yu-12-7-4-b}) that
     for $\varphi=(\varphi_1,\varphi_2)^\top\in D(\mathcal{A})$,
\begin{eqnarray*}\label{yu-12-7-11-b}
    \langle \mathcal{A}\varphi,\varphi\rangle_{\mathcal{H}}
    &=&\left\langle\left(
                   \begin{array}{c}
                     \mathcal{A}_1\varphi \\
                     \mathcal{A}_2\varphi \\
                   \end{array}
                 \right), \left(
                            \begin{array}{c}
                              \varphi_1 \\
                              \varphi_2 \\
                            \end{array}
                          \right)\right\rangle_{\mathcal{H}}\nonumber\\
                          &=&\langle (\triangle-I)\varphi_1,\varphi_1\rangle_{L^2(\mathbb{R}^N)}
                          +\langle \varphi_2(\cdot,1),\varphi_1(\cdot)\rangle_{L^2(\mathbb{R}^N)}-
                          \langle\partial_\rho\varphi_2,\varphi_2\rangle_{L^2(0,1;{L^2(\mathbb{R}^N)})}
                          \nonumber\\
                                                   &\leq&-\frac{1}{2}\|\varphi_1\|^2_{L^2(\mathbb{R}^N)}
                          +\langle \varphi_2(\cdot,1),\varphi_1(\cdot)\rangle_{L^2(\mathbb{R}^N)}
                          -\frac{1}{2}
                          \|\varphi_2(\cdot,1)\|_{L^2(\mathbb{R}^N)}^2\leq 0,
\end{eqnarray*}
    which leads to  (\ref{yu-12-7-10-b-b}).
\par
    Next, we show that  for any $\lambda>0$
        and $f\in \mathcal{H}$,
    there is a unique  $\varphi=(\varphi_1,\varphi_2)^{\top}\in D(\mathcal{A})$ such that
\begin{equation}\label{yu-12-7-12-b}
    (\lambda I-\mathcal{A})\varphi=f.
\end{equation}
    To this end, we arbitrarily fix  $\lambda>0$ and $f=(f_1,f_2)^{\top}\in \mathcal{H}$.
   Then by (\ref{yu-12-7-4-b}), we see that
    (\ref{yu-12-7-12-b}) is equivalent to the system:
\begin{equation}\label{yu-12-7-13-b}
\begin{cases}
    (\lambda I -(\triangle-I))\varphi_1(x)-\varphi_2(x,1)=f_1(x),&x\in\mathbb{R}^N,\\
    \lambda\varphi_2(x,\rho)+\tau^{-1}\partial_\rho\varphi_2(x,\rho)=f_2(x,\rho),
    &(x,\rho)\in \mathbb{R}^N
    \times(0,1),\\
    \varphi_2(x,0)=\varphi_1(x),&x\in \mathbb{R}^N.
\end{cases}
\end{equation}
  Meanwhile, after some simple computations, one can easily see that the (\ref{yu-12-7-13-b}) is equivalent to  the following system:
\begin{equation}\label{yu-12-7-14-b-b}
\begin{cases}
    \varphi_2(x,\rho)=e^{-\tau\lambda\rho}\varphi_1(x)+\tau\int_0^\rho
    e^{-\tau\lambda(\rho-s)}f_2(x,s)ds,\;\;\;\;(x,\rho)\in\mathbb{R}^N
    \times(0,1),\\
    ((\lambda+1 - e^{-\tau\lambda})I-\triangle)\varphi_1(x)=f_1(x)+\tau\int_0^1
    e^{-\tau\lambda(1-s)}f_2(x,s)ds,\;\;\;\;x\in\mathbb{R}^N.
\end{cases}
\end{equation}
   However, since
     $\lambda +1- e^{-\tau\lambda}>0$, we can use the Fourier transform to obtain that  $(\lambda+1 - e^{-\tau\lambda})I-\triangle$
    is invertible from $H^2(\mathbb{R}^N)$ to $L^2(\mathbb{R}^N)$. Therefore, the second equation
    in (\ref{yu-12-7-14-b-b})
    has a unique solution $\varphi_1$ in $H^2(\mathbb{R}^N)$.
    This, along with the first equation in (\ref{yu-12-7-14-b-b}), implies that
 (\ref{yu-12-7-12-b}) has a unique solution in $D(\mathcal{A})$.

\par
    Finally, by \eqref{yu-12-7-10-b-b} and \eqref{yu-12-7-12-b}, we can apply the Lumer-Phillips theorem (see \cite[Theorem 4.3, Chapter 1]{Pazy})
    to conclude that $\mathcal{A}$  generates a $C_0$-semigroup of contractions on
   $\mathcal{H}$. This completes the proof of  Lemma \ref{yu-lemma-12-8-1}.
\end{proof}

  Now we consider the following  system in $\mathcal{H}$:
\begin{equation}\label{yu-bbb-12-8-1}
\begin{cases}
    Y_t(t)=\mathcal{A}Y(t)+\mathcal{B}u(t)&t\in\mathbb{R}^+,\\
    Y(0)=Y_0,
\end{cases}
\end{equation}
    where $u\in L^2(\mathbb{R}^+;\mathcal{U})$ and $Y_0\in \mathcal{H}$.

\begin{remark}\label{yu-remark-12-20-2}
    By  Lemma \ref{yu-lemma-12-8-1},  we have the following two facts: Fact 1. For  each $Y_0\in\mathcal{H}$ and each
     $u\in L^2(\mathbb{R}^+;\mathcal{U})$,
      system (\ref{yu-bbb-12-8-1}) has a unique solution $Y(\cdot;Y_0,u)$ in $C([0,+\infty);
      \mathcal{H})$; Fact 2.
    The pair  $(\mathcal{A},\mathcal{B})$ satisfies  assumptions $(A_1)$-$(A_3)$ (with $(H,U)=(\mathcal{H},\mathcal{U})$).
\end{remark}
   We now give   relationships between (\ref{yu-12-7-1}) and (\ref{yu-bbb-12-8-1}), which show that (\ref{yu-12-7-1})
  can be embedded into our framework (\ref{yu-bbb-12-8-1}).

\begin{proposition}\label{yu-proposition-12-20-2}
   $(i)$
 The relationship of solutions to (\ref{yu-12-7-1}) and (\ref{yu-bbb-12-8-1}) is as follows:
  \begin{enumerate}
    \item [$(a)$] If $y(\cdot;y_0,f,u)$, with $y_0\in L^2(\mathbb{R}^N)$,
    $f\in L^2(-\tau,0;L^2(\mathbb{R}^N))$ and $u\in L^2(\mathbb{R}^+;L^2(\mathbb{R}^N))$,
         is the solution of (\ref{yu-12-7-1}), then
                            $Y(\cdot):=(Y_1(\cdot),Y_2(\cdot))^{\top}$ is the solution of
  (\ref{yu-bbb-12-8-1}) with $Y_0:=(y_0,f(-\cdot\tau))^{\top}$, where $[Y_1(t)](x):=y(t,x;y_0,f,u)$
  and $[Y_2(t)](x,\rho):=y(t-\rho\tau,x;y_0,f,u)$ for $(t,x,\rho)\in\mathbb{R}^+\times\mathbb{R}^N
  \times(0,1)$.

    \item [$(b)$] If $Y(\cdot;Y_0,u)=(Y_1(\cdot;Y_0,u),Y_2(\cdot;Y_0,u))^\top$, with
     $Y_0:=(h,g)^\top\in \mathcal{H}$ and $u\in L^2(\mathbb{R}^+;\mathcal{U})$,
          is the solution of (\ref{yu-bbb-12-8-1}), then the function $y(\cdot,\cdot)$ defined by
\begin{equation}\label{yu-12-20-8}
    y(t,x):=
\begin{cases}
    [Y_1(t;Y_0,u)](x),&(t,x)\in\mathbb{R}^+\times\mathbb{R}^N,\\
    g(x,-t\tau^{-1}),&(t,x)\in(-\tau,0)\times\mathbb{R}^N,
\end{cases}
\end{equation}
   is the solution of (\ref{yu-12-7-1}) with $y_0=h$ and $f(t,x)= g(x,-t\tau^{-1})$ for $(t,x)\in (-\tau,0)\times\mathbb{R}^N$. Moreover,
\begin{equation}\label{yu-12-20-9}
    [Y_2(t;Y_0,u)](x,\rho)=y(t-\rho\tau,x)\;\;\mbox{for}\;\;
   (t,x,\rho)\in \mathbb{R}^+\times\mathbb{R}^N\times (0,1).
\end{equation}
  \end{enumerate}

  $(ii)$ The relationship of the stabilizability of  (\ref{yu-12-7-1}) and (\ref{yu-bbb-12-8-1}) is as follows:
   System (\ref{yu-12-7-1}) is stabilizable in the sense of Definition \ref{yu-def-12-18-1} if and only if system (\ref{yu-bbb-12-8-1}) is stabilizable, i.e., there is a $K\in\mathcal{L}(\mathcal{H};\mathcal{U})$ such that $e^{(\mathcal{A}+\mathcal{B}K)t}$ is exponentially stable.
\end{proposition}
\begin{proof}
   We first prove $(a)$ of $(i)$.
   For this purpose, we arbitrarily fix  $y_0\in L^2(\mathbb{R}^N)$,
    $f\in L^2(-\tau,0;L^2(\mathbb{R}^N))$ and $u\in L^2(\mathbb{R}^+;L^2(\mathbb{R}^N))$.
    Let
\begin{equation*}\label{yu-12-18-00-1}
    z(t,x,\rho;y_0,f,u):=y(t-\rho\tau,x;y_0,f,u)\;\;\mbox{for}\;\;(t,x,\rho)\in\mathbb{R}^+
    \times\mathbb{R}^N\times(0,1).
\end{equation*}
    Then, by (\ref{yu-12-7-1}),  $(y(\cdot,\cdot;y_0,f,u),z(\cdot,\cdot,\cdot;y_0,f,u))$ satisfies  the following equation:
\begin{equation}\label{yu-12-7-3}
\begin{cases}
    y_t(t,x)=(\triangle-I)y(t,x)+z(t,x,1)+\chi_{\omega}u(t,x),&(t,x)\in\mathbb{R}^+
    \times\mathbb{R}^N,\\
    z_t(t,x,\rho)=-\tau^{-1}\partial_\rho z(t,x,\rho),&(t,x,\rho)\in\mathbb{R}^+\times\mathbb{R}^N\times(0,1),\\
    z(t,x,0)=y(t,x),&(t,x)\in\mathbb{R}^+\times\mathbb{R}^N,\\
    y(0,x)=y_0(x),&x\in\mathbb{R}^N,\\
    z(0,x,\rho)=f(-\rho\tau,x),&(x,\rho)\in \mathbb{R}^N\times(0,1).
\end{cases}
\end{equation}
    Let
    \begin{equation*}
    [Y_1(t)](x):=y(t,x;y_0,f,u);\; [Y_2(t)](x,\rho):=z(t,x,\rho;y_0,f,u)\;\;\mbox{for}\;\;(t,x,\rho)\in\mathbb{R}^+\times\mathbb{R}^N\times(0,1).
    \end{equation*}
It  is clear that
\begin{equation*}
[Y_1(0)](x):=y_0(x);\;[Y_2(0)](x,\rho):=f(-\rho\tau,x)\;\;\mbox{for}\;\;x\in\mathbb{R}^N, \rho\in(0,1).
 \end{equation*}
    Then, by (\ref{yu-12-7-4-b}) and (\ref{yu-12-7-3}), we see that
    $Y(\cdot):=(Y_1(\cdot),Y_2(\cdot))^{\top}$ is the solution to  (\ref{yu-bbb-12-8-1}), with $Y_0:=(y_0,f(-\cdot\tau))^{\top}$, which leads to  $(a)$ of $(i)$.
\par
    We next prove $(b)$ of $(i)$. To this end, we arbitrarily fix  $Y_0=(h,g)^\top\in\mathcal{H}$ and $u\in \mathcal{U}$. Let $Y(\cdot;Y_0,u)=(Y_1(\cdot;Y_0,u),Y_2(\cdot;Y_0,u))^\top$ be the solution of (\ref{yu-bbb-12-8-1}), with $Y_0=(h,g)^\top\in\mathcal{H}$.
    Then we have
    \begin{equation*}
    [Y_1(0;Y_0,u)](x)=h(x);\;[Y_2(0;Y_0,u)](x,\rho)=g(x,\rho)\;\;\mbox{for}\;\;x\in \mathbb{R}^N, \rho\in (0,1).
    \end{equation*}
    Let
    $y(\cdot,\cdot)$ be defined by (\ref{yu-12-20-8}). Then it is clear that $y(t,x):=[Y_1(t;Y_0,u)](x)$ in $(t,x)\in\mathbb{R}^+
    \times\mathbb{R}^N$. Let $z(t,x,\rho):=[Y_2(t;Y_0,u)](x,\rho)$ for $(t,x,\rho)\in\mathbb{R}^+\times
    \mathbb{R}^N\times (0,1)$.   By (\ref{yu-12-7-4-b}), we have that the pair $(y(\cdot,\cdot),z(\cdot,\cdot,\cdot))$ satisfies that equation (\ref{yu-12-7-3})
    with $y_0=h$ and $f(t,x)=g(x,-t\tau^{-1})$ for $(t,x)\in (-\tau,0)\times\mathbb{R}^N$.
    By the second, third and fifth equations in (\ref{yu-12-7-3}), we can directly check that
\begin{equation}\label{yu-12-20-10}
    z(t,x,\rho)=
\begin{cases}
    y(t-\rho\tau,x)&\mbox{if}\;\;t\geq \rho\tau,\\
   g(x,-(t-\rho\tau)\tau^{-1})&\mbox{if}\;\;t<\rho\tau,
\end{cases}
\;\;\;\;
    (t, x, \rho)\in\mathbb{R}^+\times \mathbb{R}^N\times (0,1).
\end{equation}
    This, along with the first equation in (\ref{yu-12-7-3}) and (\ref{yu-12-20-8}), implies
\begin{equation*}\label{yu-12-20-11}
    y_t(t,x)=
    (\triangle-I)y(t,x)+y(t-\tau,x)+\chi_{\omega}u(t,x), \;\;(t,x)\in\mathbb{R}^+
    \times\mathbb{R}^N.
\end{equation*}
    Since  $y(0)=Y_1(0;Y_0,u)=h$, the above shows  that the function $y(\cdot,\cdot)$, defined by (\ref{yu-12-20-8}), is the solution
    of (\ref{yu-12-7-1}), with $y_0(x)=h(x)$ and $f(t,x)= g(x,-t\tau^{-1})$ for $(t,x)\in (-\tau,0)\times\mathbb{R}^N$.
     Meanwhile,  (\ref{yu-12-20-9}) follows from (\ref{yu-12-20-8}), (\ref{yu-12-20-10}) and the definition of $z(\cdot,\cdot,\cdot)$ at once. These lead to $(b)$ of  $(i)$.
\par
     We finally prove $(ii)$.
      First, we suppose that system (\ref{yu-bbb-12-8-1}) is stabilizable, i.e., there is $\alpha>0$, $C>0$, and  $K\in\mathcal{L}(\mathcal{H};\mathcal{U})$ such that for any $Y_0\in\mathcal{H}$,
   the solution $Y^K(\cdot;Y_0)$ of the following equation:
\begin{equation*}\label{yu-12-20-30}
\begin{cases}
    Y_t(t)=\mathcal{A}Y(t)+\mathcal{B}KY(t),&t\in\mathbb{R}^+,\\
    Y(0)=Y_0,
\end{cases}
\end{equation*}
    satisfies
    \begin{equation}\label{4.28w-12-22}
    \|Y^K(t;Y_0)\|_{\mathcal{H}}\leq Ce^{-\alpha t}\|Y_0\|_{\mathcal{H}}\;\;\mbox{for any}\;\;t\geq 0.
    \end{equation}
            We arbitrarily fix  $(y_0,f)^\top\in L^2(\mathbb{R}^N)\times
    L^2(-\tau,0;L^2(\mathbb{R}^N))$. Then we define
    \begin{equation}\label{4.29-12-22-wang}
    h(x):=y_0(x);\; g(x,\rho):=f(-\rho\tau,x)\;\;\mbox{for}\;\;(x,\rho)\in\mathbb{R}^N\times(0,1).
    \end{equation}
    Let $Y_0:=(h,g)^\top$. We write
     \begin{equation}\label{4.29w-12-22}
     Y^K(t;Y_0):=(Y_1^k(t),Y_2^K(t))^\top;\;\;u^K(t):=KY^K(t;Y_0)\;\;\mbox{for}\;\;t\in\mathbb{R}^+.
      \end{equation}
     Since $K\in\mathcal{L}(\mathcal{H};\mathcal{U})$, it follows from \eqref{4.29w-12-22} and \eqref{4.28w-12-22}
     that
      $u^K\in L^2(\mathbb{R}^+;\mathcal{U})$.
      Then, according to  $(b)$ in $(i)$ of Proposition \ref{yu-proposition-12-20-2}, the function, defined by
\begin{equation}\label{yu-12-20-31}
    y^K(t,x):=
\begin{cases}
    [Y_1^K(t)](x)&(t,x)\in\mathbb{R}^+\times\mathbb{R}^N,\\
    g(x,-t\tau^{-1})&(t,x)\in(-\tau,0)\times\mathbb{R}^N,
\end{cases}
\end{equation}
    satisfies
\begin{equation*}\label{yu-12-20-32}
     y_t^K(t,x)=(\triangle-I) y^K(t,x)+y^K(t-\tau,x)+\chi_{\omega}u^K(t,x),\;\;(t,x)\in\mathbb{R}^+\times\mathbb{R}^N;
\end{equation*}
\begin{equation*}
y^K(0,x)=h(x),\;\;x\in \mathbb{R}^N;
\end{equation*}
\begin{equation*}
y^K(t,x)=g(x,-t\tau^{-1}),\;\;(t,x)\in (-\tau,0)\times\mathbb{R}^N;
\end{equation*}
\begin{equation}\label{yu-12-20-32-bb}
    [Y^K_2(t)](x,\rho)=y^K(t-\rho\tau,x)\;\;\mbox{for}\;\;
   (t,x,\rho)\in \mathbb{R}^+\times\mathbb{R}^N\times (0,1).
\end{equation}
    These, along with (\ref{yu-12-20-31}), \eqref{4.28w-12-22}
    and
    \eqref{4.29-12-22-wang}, yield
\begin{eqnarray}\label{yu-12-20-33}
    \|y^K(t)\|_{L^2(\mathbb{R}^N)}^2
    &\leq& \|Y^K(t;Y_0)\|_{\mathcal{H}}^2
    \leq C^2 e^{-2\alpha t}\left(\|h\|^2_{L^2(\mathbb{R}^N)}
    +\int_0^1\|g(\cdot,\rho)\|^2_{L^2(\mathbb{R}^N)}d\rho\right)\nonumber\\
    &\leq&C^2 e^{-2\alpha t}\left(\|h\|^2_{L^2(\mathbb{R}^N)}
    +\int_{-\tau}^0\|f(s,\cdot)\|^2_{L^2(\mathbb{R}^N)}ds\right)
    \;\;\mbox{for}\;\;t\in\mathbb{R}^+.
\end{eqnarray}

    Next, we will rewrite $u^K$ in the feedback form required by Definition \ref{yu-def-12-18-1}.
    To this end, we define an operator $\Lambda_{\tau}$ in the manner: for each $\hat f\in L^2(0,1;L^2(\mathbb{R}^N))$, we let
    \begin{equation*}
    \Lambda_{\tau}[\hat f](t,x):=\hat f(x,-t\tau^{-1}),\;\;(t,x)\in(-\tau,0)\times\mathbb{R}^N.
     \end{equation*}
     One can directly check that
    $\Lambda_{\tau}: L^2(0,1;L^2(\mathbb{R}^N))\to L^2(-\tau,0;L^2(\mathbb{R}^N))$
     is isomorphic.
     Then, by \eqref{4.29w-12-22} and  (\ref{yu-12-20-32-bb}), we find that
\begin{equation*}\label{yu-12-20-34}
    u^K(t)=K\left(
              \begin{array}{cc}
                I & 0 \\
                0 & \Lambda_{\tau}^{-1} \\
              \end{array}
            \right)\left(
                     \begin{array}{c}
                       y^K(t) \\
                       y^K(t+\cdot) \\
                     \end{array}
                   \right).
\end{equation*}
    The above, along with  (\ref{yu-12-20-33}) and  Definition \ref{yu-def-12-18-1}, shows that
            system (\ref{yu-12-7-1}) is stabilizable in the sense of Definition \ref{yu-def-12-18-1}
    with the feedback law $K\left(
              \begin{array}{cc}
                I & 0 \\
                0 & \Lambda^{-1} \\
              \end{array}
            \right)$.

\par
   Conversely, we suppose that system (\ref{yu-12-7-1}) is stabilizable in the sense of Definition \ref{yu-def-12-18-1}, i.e.,  there is $\alpha>0$, $C>0$ and a feedback law $K\in\mathcal{L}(L^2(\mathbb{R}^N)
    \times L^2(-\tau,0;L^2(\mathbb{R}^N));L^2(\mathbb{R}^N))$ such that  for any $y_0\in L^2(\mathbb{R}^N)$ and $f\in L^2(-\tau,0;L^2(\mathbb{R}^N))$, the solution $y_K(\cdot;y_0,f)$ of the closed-loop system (\ref{yu-12-20-100}) satisfies (\ref{yu-12-20-101}).
   We arbitrarily fix  $Y_0=(h,g)^\top\in\mathcal{H}$. Let
    \begin{equation*}\label{4.34-ww-12-22}
    y_0(x)=h(x);\;\;f(t,x)=g(x,-t\tau^{-1}),\;\;(x,t)\in\mathbb{R}^N\times(-\tau,0).
    \end{equation*}
    Then
     $(y_0,f)\in L^2(\mathbb{R}^N)\times L^2(-\tau,0;L^2(\mathbb{R}^N))$. Let
      \begin{equation}\label{4.34-12-22-w}
      u_K(t):=K(y_K(t;y_0,f),y(t+\cdot;y_0,f))^\top,\;\;t\geq 0.
      \end{equation}
      By \eqref{4.34-12-22-w}, one can easily check that
      \begin{equation}\label{4.35-12-22-w}
    y_K(t,x;y_0,f)=y(t,x;y_0,f,u_K),\;\;t\geq 0, x\in \mathbb{R}^N.
\end{equation}
   Given $t\geq 0$, $x\in \mathbb{R}^N$ and $\rho\in(0,1)$, we let
    \begin{equation}\label{4.37-ww-12-22}
    [Y^K_1(t)](x):=y(t,x;y_0,f,u_K);\;\;[Y^K_2(t)](x,\rho):=y(t-\rho\tau,x;y_0,f,u_K).
    \end{equation}
     By \eqref{4.37-ww-12-22}, \eqref{4.35-12-22-w} and  $(a)$ in $(i)$ of Proposition \ref{yu-proposition-12-20-2}, we see that $Y^K(\cdot):=(Y^K_1(\cdot),Y^K_2(\cdot))^{\top}$ is the solution of
  (\ref{yu-bbb-12-8-1}), with $Y_0:=(y_0,f(-\cdot\tau))^{\top}$ and $u=u^K$.
  This, together with \eqref{4.34-12-22-w}, \eqref{4.37-ww-12-22} and
   (\ref{yu-12-20-101}), yields
\begin{equation*}\label{yu-12-20-36}
    \|Y^K(t)\|^2_{\mathcal{H}}=\|y_K(t)\|^2_{L^2(\mathbb{R}^N)}
    +\int_0^1\|y_K(t-\rho\tau,\cdot;y_0,f)\|^2_{L^2(\mathbb{R}^N)}d\rho
    \leq C^2(1+e^{2\tau})e^{-2\alpha t}\|Y_0\|^2_{\mathcal{H}}\;\;\mbox{for}\;\;t\in\mathbb{R}^+.
\end{equation*}
    With \eqref{4.34-12-22-w}, the above leads to
    \begin{equation}\label{4.38-ww-12-22}
    \mathcal{U}_{ad}(Y_0):=\{u\in L^2(\mathbb{R}^+;\mathcal{U}): Y(\cdot;Y_0,u)\in L^2(\mathbb{R}^+;\mathcal{H})\}\neq \emptyset,
    \end{equation}
         where $Y(\cdot;Y_0,u)$ is the solution of (\ref{yu-bbb-12-8-1}), with the initial value $Y_0$ and the control $u$.
    By \eqref{4.38-ww-12-22}, we can apply   \cite[Proposition 3.9]{Liu-Wang-Xu-Yu} to see that  system  (\ref{yu-bbb-12-8-1})
    is stabilizable.
\par
    Thus, we complete the proof of Proposition \ref{yu-proposition-12-20-2}.
\end{proof}

 Next, we will show that  system (\ref{yu-12-7-1}), with $u=0$, is not exponentially stable in $L^2(\mathbb{R}^N)$.
This shows the importance of studying the stabilizability of system (\ref{yu-12-7-1}).

\begin{proposition}\label{yu-proposition-12-15-1}
   System
     (\ref{yu-12-7-1}), with $u=0$, is not exponentially stable.
\end{proposition}
\begin{proof}
First of all, according to $(i)$ in Proposition \ref{yu-proposition-12-20-2}, system (\ref{yu-12-7-1}) (with $u=0$) is exponentially stable
if and only if system (\ref{yu-bbb-12-8-1}) (with $u=0$) is exponentially stable.

Now, we let $\lambda_j:=j^{-1}$ ($j\in\mathbb{N}^+$). Then let   $f^j:=(f^j_1,0)\in \mathcal{H}$,
    with $\|f^j\|_{\mathcal{H}}=1$ and $\mbox{supp}(\mathcal{F}[f_1^j])
    \subset \{\xi\in\mathbb{R}^N||\xi|\leq \sqrt{\lambda_j}\}$. It is obvious that
    $\lambda_j+1-e^{-\tau\lambda_j}>\lambda_j>0$
    for all $j\in\mathbb{N}^+$.
    Thus, equation  (\ref{yu-12-7-13-b}) (equivalently, (\ref{yu-12-7-14-b-b})),
     with $\lambda=\lambda_j$ and $f=f^j$,
    has a unique solution $\varphi^j:=(\varphi^j_1,\varphi^j_2)$ in $D(\mathcal{A})$.
    Applying the Fourier
    transform to the second equation in (\ref{yu-12-7-14-b-b}) (with $\lambda=\lambda_j$),
    using the fact that $f_2=0$ and the Plancherel theorem, we obtain that as $j\to+\infty$,
\begin{eqnarray*}
    \|\varphi^j\|_{\mathcal{H}}^2&\geq& \|\varphi_1^j\|^2_{L^2(\mathbb{R}^N)}
    =\int_{\mathbb{R}^N}
    (\lambda_j+1-e^{-\tau\lambda_j}+|\xi|^2)^{-2}|\mathcal{F}[f_1^j](\xi)|^2d\xi\nonumber\\
    &\geq& (2\lambda_j+1-e^{-\tau\lambda_j})^{-2}
    \int_{\mathbb{R}^N}|\mathcal{F}[f_1^j](\xi)|^2d\xi=(2j^{-1}+(1-e^{-\tau j^{-1}}))^{-1}\to+\infty.
\end{eqnarray*}
    This yields that $\|(\lambda_j I-\mathcal{A})^{-1}
   \|_{\mathcal{L}(\mathcal{H})}
    \to +\infty$ as $j\to +\infty$. Thus, we have
    \begin{eqnarray*}
    \sup_{\lambda\in \mathbb{C},\mbox{Re}\lambda>0}
    \|(\lambda I-\mathcal{A})^{-1}\|_{\mathcal{L}(\mathcal{H})
    }=+\infty.
    \end{eqnarray*}
    The above, along with
     \cite[Theorem 1.11, Chapter V]{Engel-Nagel}, yields  that
       (\ref{yu-bbb-12-8-1}) (with $u=0$) is not exponentially stable.
    This completes the proof.
\end{proof}

\par
    The main result of this subsection concerns with the stabilizability of system (\ref{yu-12-7-1}), which is proved
   by making use of Theorem
    \ref{yu-theorem-5-20-1}.

\begin{theorem}\label{yu-theorem-12-8-1}
    If  $\omega$ is a thick set, then system (\ref{yu-12-7-1}) is  stabilizable in the sense of Definition \ref{yu-def-12-18-1}.
\end{theorem}
\begin{proof}
    According to \cite[Theorem 1.1 and Remark $(b_1)$]{Wang-Zhang},
          the control heat equation
     $y_t=\triangle y+\chi_{\omega}u$ (in
     $\mathbb{R}^+\times \mathbb{R}^N$), with  $\omega$  a thick set, is null controllable, and consequently is rapidly stabilizable.
      Thus,
     by $(ii)$ in Theorem \ref{yu-theorem-7-3-1},
     we can conclude that for each $\gamma>0$, there is $C_0(\gamma)>0$ such that
\begin{equation}\label{yu-12-9-1-b}
    \|\varphi\|^2_{L^2(\mathbb{R}^n)}\leq C_0(\gamma)
    (\|(\lambda I -\triangle)\varphi\|^2_{L^2(\mathbb{R}^N)}
    +\|\chi_{\omega}\varphi\|_{L^2(\mathbb{R}^N)}^2)
    \;\;\mbox{for}\;\;\lambda\in \mathbb{C}^+_{-\gamma},\;\;
    \varphi\in H^2(\mathbb{R}^N).
\end{equation}
Meanwhile, one can easily check that $\mathcal{B}^*=(\chi_{\omega},0)$ and
\begin{equation}\label{yu-12-7-7}
\begin{cases}
    \mathcal{A}^*\varphi:=\left(
                          \begin{array}{c}
                            \mathcal{A}^{\sharp}_1\varphi \\
                            \mathcal{A}^{\sharp}_2\varphi \\
                          \end{array}
                        \right)\;\;\mbox{for}\;\;\varphi=\left(
                                                  \begin{array}{c}
                                                    \varphi_1 \\
                                                    \varphi_2 \\
                                                  \end{array}
                                                \right)
                      \in D(\mathcal{A}^*),\\
                        \;\;\mbox{with}\;\;\left(
                          \begin{array}{c}
                            [\mathcal{A}^{\sharp}_1\varphi](x) \\
                            [\mathcal{A}^{\sharp}_2\varphi](x,\rho) \\
                          \end{array}
                        \right)=
               \left(
                         \begin{array}{c}
                           (\triangle-I)\varphi_1(x)+\varphi_2(x,0) \\
                            -\tau^{-1}\partial_\rho\varphi_2(x,\rho)\\
                         \end{array}
                       \right)\;\;(x,\rho)\in\mathbb{R}^N\times(0,1),\\
                       D(\mathcal{A}^*):=
    \{(f_1,f_2)^{\top}\in \mathcal{H}: f_1\in  H^2(\mathbb{R}^N),
    \;f_2\in H^1(0,1;L^2(\mathbb{R}^N)),
    \;f_1(x)=f_2(x,1)\;\mbox{a.e.}\;x\in\mathbb{R}^N\}.
\end{cases}
\end{equation}

We arbitrarily fix $\gamma_0\in \mathbb{R}^+$. We claim that the pair $(\mathcal{A}^*,\mathcal{B}^*)$ satisfies (HESI)$_{\gamma_0}$, i.e.,
    there is $C(\gamma_0)>0$ such that
\begin{equation}\label{yu-12-8-4}
    \|\varphi\|^2_{\mathcal{H}}
    \leq C(\gamma_0)(\|(\lambda I-\mathcal{A}^*)\varphi\|^2_{\mathcal{H}}
    +D(\gamma_0)\|\mathcal{B}^*\varphi\|_{\mathcal{U}}^2)\;\;\mbox{for}\;\;
    \lambda\in\mathbb{C}_{-\gamma_0}^+,\;\;
    \varphi\in D(\mathcal{A}^*).
\end{equation}
    For this purpose, we arbitrarily fix
    $\lambda\in\mathbb{C}^+_{-\gamma_0}$ and $\varphi=(\varphi_1,\varphi_2)^\top \in
    D(\mathcal{A}^*)$. Let
    \begin{equation}\label{4.42-12-23}
    f=(f_1,f_2)^{\top}:=(\lambda I-\mathcal{A}^*)\varphi.
    \end{equation}
         By (\ref{yu-12-7-7}), we have
\begin{equation}\label{yu-12-8-5}
\begin{cases}
    ((\lambda+1)I-\triangle)\varphi_1(x)-\varphi_2(x,0)=f_1(x),&x\in
    \mathbb{R}^N,\\
    (\lambda I -\tau^{-1}\partial_\rho)\varphi_2(x,\rho)=f_2(x,\rho),&(x,\rho)\in
    \mathbb{R}^N\times(0,1),\\
    \varphi_1(x)=\varphi_2(x,1),&x\in\mathbb{R}^N.
\end{cases}
\end{equation}
    From the second and third equations in (\ref{yu-12-8-5}), we deduce that for each $\rho\in(0,1)$,
\begin{equation}\label{yu-12-8-6}
    \varphi_2(x,\rho)=e^{-\tau\lambda(1-\rho)}\varphi_1(x)
    +\tau\int_\rho^1e^{-\tau\lambda(s-\rho)}f_2(x,s)ds,\;\;\;\;x\in\mathbb{R}^N.
\end{equation}
    The combination of \eqref{yu-12-8-6} and the first equation in (\ref{yu-12-8-5}) leads to
\begin{equation}\label{yu-12-8-7}
    ((\lambda+1-e^{-\tau\lambda})I-\triangle)\varphi_1(x)=f_1(x)
+\tau\int_0^1
    e^{-\tau\lambda s}f_2(x,s)ds,\;\;\;\;x\in\mathbb{R}^N.
\end{equation}
   Since
   \begin{equation*}
   \mbox{Re}(\lambda+1-e^{-\tau\lambda})\geq \mbox{Re}\lambda
    +1-e^{-\tau\mbox{Re}\lambda}\geq -(\gamma_0
    +e^{\tau\gamma_0}-1)\;\mbox{for each}\;\lambda\in \mathbb{C}_{-\gamma_0}^+
   \end{equation*}
   and
   \begin{equation*}
   \gamma_0
    +e^{\tau\gamma_0}-1>0,
   \end{equation*}
   the combination of (\ref{yu-12-9-1-b}) (with $\gamma:=\gamma_0
    +e^{\tau\gamma_0}-1$) and (\ref{yu-12-8-7}) yields
     that there is $C_1(\gamma_0)>0$ such that
\begin{eqnarray}\label{yu-12-8-8}
    \|\varphi_1(\cdot)\|_{L^2(\mathbb{R}^N)}^2&\leq& C_1(\gamma_0)
    \left(\left\|f_1(\cdot)
    +\tau\int_0^1
    e^{-\tau\lambda s}f_2(\cdot,s)ds\right\|_{L^2(\mathbb{R}^N)}^2
    +\|\chi_{\omega}\varphi_1(\cdot)\|^2_{L^2(\mathbb{R}^N)}\right)\nonumber\\
    &\leq& 2C_1(\gamma_0)(1+\tau  e^{2\tau\gamma_0})\|f\|_{\mathcal{H}}^2
    +C_1(\gamma_0)\|\mathcal{B}^*\varphi\|_{\mathcal{U}}^2.
\end{eqnarray}
   Meanwhile, it follows from (\ref{yu-12-8-6}) that
\begin{eqnarray}\label{yu-12-8-11}
    \|\varphi_2(\cdot,\cdot)\|^2_{L^2(0,1;L^2(\mathbb{R}^N))}&\leq&
    2e^{2\tau\gamma_0}\|\varphi_1(\cdot)\|_{L^2(\mathbb{R}^N)}^2
    +2\tau^2 e^{2\tau\gamma_0}\int_0^1\|f_2(\cdot,s)\|_{L^2(\mathbb{R}^N)}^2ds\nonumber\\
    &\leq& \left[4 C_1(\gamma_0)e^{2\tau\gamma_0}(1+\tau e^{2\tau\gamma_0})
    +2\tau e^{2\tau\gamma_0}\right]\|f\|_{\mathcal{H}}^2\nonumber\\
    &\;&+2 C_1(\gamma_0)e^{2\tau\gamma_0}
    \|\mathcal{B}^*\varphi\|_{\mathcal{U}}^2.
\end{eqnarray}
   Now, by \eqref{4.42-12-23},
     (\ref{yu-12-8-8})
    and (\ref{yu-12-8-11}), we obtain (\ref{yu-12-8-4}), with
    \begin{eqnarray*}
    C(\gamma_0):=4C_1(\gamma_0)(1+\tau e^{2\tau\gamma_0})
    (1+\tau e^{2\tau\gamma_0})
    +2\tau^2e^{2\tau\gamma_0},\;\;  D(\gamma_0):=C_1(\gamma_0)(1+2e^{2\tau\gamma_0}).
    \end{eqnarray*}

\par
  Finally, according to Lemma \ref{yu-lemma-12-8-1},
    the $C_0$-semigroup generated by $\mathcal{A}^*$ is uniformly bounded, which leads to the condition
    $(a)$ in Theorem \ref{yu-theorem-5-20-1}. Thus, we can apply Theorem
    \ref{yu-theorem-5-20-1} and (\ref{yu-12-8-4}) to conclude
         that  system \eqref{yu-bbb-12-8-1} is stabilizable. Then by $(ii)$ of Proposition \ref{yu-proposition-12-20-2}, we find  that  system (\ref{yu-12-7-1}) is stablilizable in the sense of Definition \ref{yu-def-12-18-1}.
   This completes the proof.
\end{proof}

\section{Appendix}
 \begin{proposition}\label{yu-prop-6-3-1}
    Suppose that $(A_1)$-$(A_3)$ hold. Then the inequalities  \eqref{yu-6-9-1} and \eqref{yu-6-12-bb-0} are equivalent.
\end{proposition}
\begin{proof}
    We first show \eqref{yu-6-9-1}$\Rightarrow$\eqref{yu-6-12-bb-0}. Suppose that
    \eqref{yu-6-9-1} holds. Let $\beta_1\in (0,\beta)$.
    Then, it follows from  (\ref{yu-6-9-1}) that for $\lambda\in\mathbb{C}^+_{-\beta_1}$,
\begin{equation*}\label{yu-6-3-1}
    \|\varphi\|_H^2\leq \frac{C(\beta)}{(\beta-\beta_1)^2}\left(\|(\lambda I-A^*)\varphi\|_H^2
    +\|B^*\varphi\|_U^2\right)\;\;\mbox{for}\;\;\varphi\in H_1,
\end{equation*}
   which leads to \eqref{yu-6-12-bb-0} with a different $C(\beta)>0$.

   Next, we  show \eqref{yu-6-12-bb-0}$\Rightarrow$\eqref{yu-6-9-1}. Suppose that \eqref{yu-6-12-bb-0} is true.
   First,
   there are two constants $\omega>0$ and $C(\omega)>0$  such that
        $\|S^*(t)\|\leq C(\omega)e^{\omega t}$ for all $t\in\mathbb{R}^+$,
        which implies that  $\mathbb{C}_{\omega}^+\subset \rho(A^*)$ and that  for each  $\lambda\in \mathbb{C}_{\omega}^+$, $\|(\lambda I-A^*)^{-1}\|_{\mathcal{L}(H)}\leq C(\omega)(\mbox{Re}\lambda-\omega)^{-1}$
        (see \cite[Theorem 5.3 and Remark 5.4, Section 1.5, Chapter 1] {Pazy}). These facts, together with the same argument in \eqref{yu-8-17-1}, imply that for each $\lambda\in \mathbb{C}_{\max\{\omega,2|\beta-\omega|-\beta\}}^+$,
\begin{eqnarray}\label{yu-6-3-2}
    \|\varphi\|_H\leq \frac{2C(\omega)}{\mbox{Re}\lambda+\beta}
    \|(\lambda I-A^*)\varphi\|_H \;\;\mbox{for}
    \;\;\varphi\in H_1.
\end{eqnarray}
   Meanwhile, it follows from \eqref{yu-6-12-bb-0} that for  $\lambda\in \mathbb{C}_{-\beta}^+\setminus\mathbb{C}_{\max\{\omega,2|
     \beta-\omega|-\beta\}}^+$,
     we have
\begin{eqnarray*}\label{yu-6-3-3}
    \|\varphi\|_H^2&\leq& \frac{(\beta+\max\{\omega,2|
     \beta-\omega|-\beta\})^2}{(\mbox{Re}\lambda+\beta)^2}C(\beta)
     \left(\|(\lambda I-A^*)\varphi\|_H^2+
     \|B^*\varphi\|_U^2\right)\;\;\mbox{for}\;\;\varphi\in H_1,
\end{eqnarray*}
    where we recall that the quotient in front of $C(\beta)$ are large than or equal to $1$. The combination of this result and (\ref{yu-6-3-2}) yield
\begin{equation*}\label{yu-6-12-3}
    \|\varphi\|_H^2\leq \frac{C(\beta,\omega)}{(\mbox{Re}\lambda+\beta)^2}
    \left(\|(\lambda I-A^*)\varphi\|_H^2+\|B^*\varphi\|_U^2\right)
    \;\;\mbox{for}\;\;\lambda\in\mathbb{C}^+_{-\beta},\;\varphi\in H_1,
\end{equation*}
    where $C(\beta,\omega):=(\beta+\max\{\omega,2|
     \beta-\omega|-\beta\})^2C(\beta)+4(C(\omega))^2$.
     This implies that \eqref{yu-6-9-1} holds.
     This completes the proof.
\end{proof}

\begin{remark}\label{yu-remark-7-31-1}
     The following example shows that \eqref{yu-6-9-1} is more sharp than \eqref{yu-6-12-bb-0} to describe the optimal decay rate of system $[A,B]$: Suppose that $A\in\mathbb{R}^{n\times n}$ and $B\in\mathbb{R}^{n\times m}$ ($n,m\in\mathbb{N}^+$), i.e., system $[A,B]$ is a finite-dimensional system in $\mathbb{R}^n$. Further, we assume that $[A,B]$ is stabilizable, but not controllable. It follows that there is an invertible matrix $P\in\mathbb{R}^{n\times n}$ such that
$$
    PAP^{-1}=\left(
               \begin{array}{cc}
                 A_1 & A_2 \\
                 0 & A_3 \\
               \end{array}
             \right),\;\;PB=\left(
                             \begin{array}{c}
                               B_1 \\
                               0 \\
                             \end{array}
                           \right),
$$
    where $[A_1,B_1]$ is controllable and $\sigma(A_3)$ is non-empty and in $\mathbb{C}_0^{-}$, i.e.,
    $\sigma^*:=\max\{\mbox{Re}\,\lambda:\lambda\in \sigma(A_3)\}<0$.
    We define the optimal decay rate of system $[A,B]$ as follows:
$$
    \sigma^{\sharp}:=\inf\{\alpha\in\mathbb{R}:\exists F\in\mathbb{R}^{m\times n}\;
    \mbox{s.t.}\;A+BF\;\mbox{is exponentially stable with the decay rate}\;\alpha\}.
$$
    One can directly check that $\sigma^{\sharp}=-\sigma^*$. If there is $F^\sharp\in\mathbb{R}^{m\times n}$ such that $A+BF^\sharp$ is exponentially stable with decay rate $\alpha^\sharp$, then we say the optimal decay rate of $[A,B]$ can be reached, otherwise, we say that it can not be  reached. We take $\lambda^*\in \sigma(A_3)(\subset \sigma(A))$ such that $\mbox{Re}\lambda^*=\sigma^*$. If the geometric multiplicity of $\lambda^*$ equals to its algebraic multiplicity, then, by the classical argument, we can directly check that the optimal decay rate of $[A,B]$ can be reached. Moreover, by the Laplace transform, we can conclude that \eqref{yu-6-9-1} holds for $\beta=\sigma^{\sharp}$, but \eqref{yu-6-12-bb-0} holds only for $\beta<\sigma^{\sharp}$.
    If the geometric multiplicity of $\lambda^*$ is strictly less than  its algebraic multiplicity, then, we can directly check that the optimal decay rate of $[A,B]$ can not be reached. In this case, we can show that  \eqref{yu-6-9-1} and \eqref{yu-6-12-bb-0} hold only for $\beta<\sigma^{\sharp}$. In summary, \eqref{yu-6-9-1} holds for $\beta=\alpha^\sharp$ in  some cases, while \eqref{yu-6-12-bb-0} holds only for $\beta<\alpha^\sharp$. Therefore, we say that \eqref{yu-6-9-1} is more sharp than \eqref{yu-6-12-bb-0}.
\end{remark}

\end{document}